\documentclass[a4paper,twoside,12pt]{amsart} 
 \usepackage[T1]{fontenc}            
 \usepackage{epsf}
 \usepackage{amscd}                  
 \usepackage{xypic}                  
 \usepackage{amssymb}
\begin{document}

\pagenumbering{arabic}
  \renewcommand{\datename}{\textit{Dato}:} 
\def\Im{\operatorname{Im}}          
\newtheorem{thm}{Theorem}[section]
\newtheorem{kor}[thm]{Korollar}
\newtheorem{cor}[thm]{Corollary}
\newtheorem{construction}[thm]{Construction}
\newtheorem{conj}[thm]{Conjecture}
\newtheorem{lemma}[thm]{Lemma}
\newtheorem{prop}[thm]{Proposition}
\newtheorem{oppg}[thm]{Oppgave}
\theoremstyle{definition}
\newtheorem{defn}[thm]{Definition}
\newtheorem{bem}[thm]{Remark}
\newtheorem{notation}[thm]{Notation}
\newtheorem{example}[thm]{Example}
\renewcommand{\proofname}{Proof}
\let\u=\underline
\newcommand{\vcgreek}[1]{\boldsymbol{#1}} 
\def\i{^{-1}}
\def\RR{{\mathbb R}}
\def\kar{\operatorname{char}}
\def\coker{\operatorname{coker}}
\def\mod{\operatorname{mod}}
\def\dim{\operatorname{dim}}
\def\ker{\operatorname{ker}}
\def\deg{\operatorname{deg}}
\def\konstant{\operatorname{konstant}}
\def\ann{\operatorname{ann}}
\def\lin{\operatorname{lin}}
\def\Sing{\operatorname{Sing}}
\def\min{\operatorname{min}}
\def\H{\operatorname{H}}
\def\depth{\operatorname{depth}}
\def\pd{\operatorname{pd}}
\def\im{\operatorname{Im}}
\def\rank{\operatorname{rank}}
\def\Sec{\operatorname{Sec}}
\def\Res{\operatorname{Res}}
\def\|{\mid}
\def\CC{{\mathbb C}}
\def\GG{{\mathbb G}}
\def\ZZ{{\mathbb Z}}
\def\DD{{\mathbb D}}
\def\NN{{\mathbb N}}
\def\QQ{{\mathbb Q}}
\def\VV{{\mathbb V}}
\def\PP{{\mathbb P}}
\def\CM{Cohen-Macaulay }
\def\FF{{\mathbb F}}
\def\AA{{\mathbb A}}
\def\D{{\mathcal D}}
\def\A{{\mathcal A}}
\def\F{{\mathcal F}}
\def\J{{\mathcal J}}
\def\G{{\mathcal G}}
\def\M{{\mathcal M}}
\def\T{{\mathcal T}}
\def\N{{\mathcal N}}
\def\O{{\mathcal O}}
\def\I{{\mathcal I}}
\def\Q{{\mathcal Q}}
\def\E{{\mathcal E}}
\def\K{{\mathcal K}}
\def\C{{\mathcal C}}
\def\Z{{\mathcal Z}}
\def\V{{\mathcal V}}
\def\B{{\mathcal B}}
\def\L{{\mathcal L}}
\def\iso{\cong}
\def\congr{\equiv}
\def\sub{\subseteq}
\def\subne{\subsetneqq}
\def\hpil{\longrightarrow}
\def\surj{\hpil\mspace{-26.0mu}\hpil}
\def\Aut{\operatorname{Aut}}
\def\id{\operatorname{id}}
\def\Der{\operatorname{Der}}
\def\Tor{\operatorname{Tor}}
\def\Ext{\operatorname{Ext}}
\def\red{\operatorname{red}}
\def\spec{\operatorname{Spec}}
\def\Proj{\operatorname{Proj}}
\def\Hom{\operatorname{Hom}}
\def\End{\operatorname{End}}
\def\Hilb{\operatorname{Hilb}}
\def\Grass{\operatorname{Grass}}
\def\Pic{\operatorname{Pic}}
\def\Supp{\operatorname{Supp}}
\def\Sym{\operatorname{Sym}}
\def\GL{\operatorname{GL}}
\def\SL{\operatorname{SL}}
\def\Pic{\operatorname{Pic}}
\def\codim{\operatorname{codim}}
\def\nil{\operatorname{nil}}
\def\dh{\operatorname{dh}}
\def\det{\operatorname{det}}
\hfuzz5pc 

\title {Multiple Structures}

\author {Jon Eivind Vatne}

\address{ Matematisk institutt\\
          Johs. Brnsgt. 12\\
          N-5008 Bergen\\
          Norway}

\email{jonev@mi.uib.no}

\begin{abstract}
We give a systematic approach to constructing non-reduced, locally
Cohen-Macaulay schemes, mostly with reduced support a projective
variety.  The hierarchy of such structures includes a lot of
information about the variety, its embedding in projective space and
the behaviour of its vector bundles.  For instance, Hartshorne's
Conjecture on complete intersections in codimension two is
reformulated in terms of existence of certain schemes of degrees two
and three.  We also give a lot of examples of multiple structures, and
lists of multiple structures on linear subspaces of codimension two in
low degree.
\end{abstract}

\maketitle

\section{\label{intro} Introduction}

The aim of this paper is present a systematic approach to the
classification theory of (locally) Cohen-Macualay nilpotent schemes
whose underlying space is a smooth, projective variety.  On the way,
we get lots of interesting examples, and several results linking these
structures to more general questions of projective geometry.\\

In Chapter \ref{s1} we present the main technical tools needed for the
theory.  Basically, for any smooth projective subvariety $X\subset
\PP^N$ and every nilpotent scheme $Y$ with $Y_{\red}=X$, we intersect
$Y$ with the various infinitesimal neighbourhoods of $X$, and then
remove the embedded components.  This gives us a filtration of $Y$,
which we will call the {\em $S_1$-filtration} of $Y$.  There are very
strong connections between the members of this filtration.  For our
purposes, the most interesting point is that $Y$ can be reconstructed
from short exact sequences between the ideal sheaves of the schemes in
the filtration and certain $\O_X$-Modules.  When the multiple schemes
are Cohen-Macaulay, these $\O_X$-Modules are locally free, and
vice-versa.\\

In Chapter \ref{codimtwo},  we give a complete list of the
Cohen-Macaulay schemes whose reduced subscheme is a linear subspace of
codimension two, up to and including multiplicity five, with the extra
assumption that all members of the $S_1$-filtration are Cohen-Macaulay
themselves.  We also get some general results, showing how the
characteristic and the dimension come into play.\\

In Chapter \ref{kvotient-->multippelstruktur} we give results linking
the classification theory of multiple structures to other questions in
projective geometry.  The main observation is that if the locally free
$\O_X$-Module metioned above is an extension of locally free sheaves
of smaller rank, there is an intermediate Cohen-Macaulay nilpotent
structure.  This allows us to reformulate Hartshorne's Conjecture on
complete intersections in codimension two entirely in terms of
existence questions for schemes of degree two and three.\\

In Chapter \ref{ikkeI} we present a large class of 2-dimensional
multiple structures in $\PP^4$ that have non-Cohen-Macaulay schemes in
the $S_1$-filtration.  They all contain a complete intersection of
multiplicity two less.\\

{\bf Acknowledgements}

This article is based on parts of my Doctorate thesis.  I would like
to thank the group in algebraic geometry at the University of Bergen,
G.Fløystad, T.Johnsen, A.L.Knutsen, S.A.Strømme, J.-M. Økland and most
improtantly my advisor A.Holme.  Furthermore, I have received a lot of
interesting comments and advice from the evalutation committee:
T.Johnsen, K.Ranestad and F.-O.Schreyer.

\section{\label{s1} Recursive construction of multiple structures}

  In this chapter we will construct a sequence of subschemes for any 
  nilpotent scheme, whose reduced subscheme is a smooth subvariety of a
  smooth variety.  The subschemes have milder nilpotencies, and
  satisfy Serre's $(S_1)$-criterion.  However, if the nilpotent scheme
  is Cohen-Macaulay, or a (locally) complete intersection, the same is
  not necessarily true for these subschemes.  Note that we always use
  Cohen-Macaulay in the EGA sense: a scheme  is
  Cohen-Macaulay if all its local rings are Cohen-Macaulay rings (and
  similarly for Modules).\\

\subsection{Definition and elementary properties of the $S_1$-filtration}
Let $X \subset Y \subset Z$, where $X\subset Z$ is a regular closed
embedding of smooth connected projective schemes, and $Y$ is a
nilpotent scheme with $Y_{\red}=X$.  $\I_X,\I_Y$ are the
Ideals, $I_X,I_Y$ the corresponding (saturated) ideals in the (or a)
homogeneous coordinate ring of the projective variety $Z$.  The
infinitesimal neighbourhoods of $X$, $X^{(i)}$, are defined by the
Ideals $\I_X^{i+1}$.  Since $Y_{\red}=X$, there is a unique
integer $k$ such that $Y$ is contained in $X^{(k)}$ but is not
contained in $X^{(k-1)}$.\\

\begin{defn}
For each integer $j\leq k$, define the scheme $Y_j$ as follows:
Intersect $Y$ with $X^{(j)}$ and remove the embedded components.
Algebraically, $I_j:=I_{Y_j}$ is the $I_X$-primary component of
$I_Y+I_X^{j+1}$.  Thus we have a sequence of subschemes   

\begin{equation}
\label{s1filtrasjonen}
X=Y_0 \subset Y_1 \subset ... \subset Y_{k-1} \subset Y_k \subset Y
\end{equation} 
which we will call {\em the $S_1$-filtration of $Y$}.
\end{defn}

  The name is not misleading:

\begin{prop}
The schemes $Y_j$ all satisfy $(S_1)$.
\end{prop}
\begin{proof}
According to EGA Proposition $IV_2,\S 5.7.5$, a scheme satisfies
Serre's $(S_1)$-criterion if and only if it is without embedded
components.  But we have removed the embedded components from our
scheme.
\end{proof}

\begin{bem}
\label{YS1}
If $Y$ itself satisfies $(S_1)$, the last inclusion $Y_k\subset Y$
actually becomes an equality, and we have a {\em filtration} of $Y$.
\end{bem}

\begin{bem}
The projectivity condition is included for simplicity.  By considering
only the Ideal sheaves, most of the theory that follows can be
generalized to the setting of smooth, connected $X$ and $Z$, not even
necessarily quasi-projective.  The filtration can also be defined
without smoothness conditions on $X$ and $Z$, but then the results must be
substantially weakened.
\end{bem}

Thus the schemes $Y_j$ are not necessarily Cohen-Macaulay, even if $Y$
is, just
locally of depth $\geq 1$.  But we have a partial result:

\begin{prop}\label{cmcod1}
The locus $D_j\subset X$ where $Y_j$ fails to be Cohen-Macaulay has codimension
at least two.
\end{prop}
\begin{proof}
Whenever $x\in X$ satisfies $\codim \overline{\{x\}}=1$, $Y_j$ is
Cohen-Macaulay in $x$ since $Y_j$ satisfies $(S_1)$.  But being
Cohen-Macaulay is an open condition, so in each irreducible subset of
codimension one there is an open dense set where $Y_j$ is
Cohen-Macaulay.  Thus $\codim D_j \geq 2$.
\end{proof}

\begin{cor}\label{s1kurve}
 If $X$ is a curve, $Y_j$ is Cohen-Macaulay.
\end{cor}

\begin{defn} A multiple Cohen-Macaulay scheme $Y$ is called a multiple
  structure of {\em type I} if the $S_1$-filtration consists of
  Cohen-Macaulay schemes.  So all multiple structures on smooth curves
  are of type I.  This definition comes from Manolache,\cite{Man92}.
\end{defn}

For our study of the terms of the filtration we need the following
technical lemma:

\begin{lemma}\label{alglemma} Suppose $R$ is a noetherian ring whose zero ideal is
  without embedded associated primes.  Let $J \subset R$ be an ideal with
  $J_P=0$ for all minimal primes $P$ in $R$.  Then $J=0$.
\end{lemma}
\begin{proof}
Assuming $J_P=0$ for all minimal primes $P$, we can find an element $s\in R$ which is not contained in any minimal
prime $P$, and which annihilates $J$ (e.g. by prime avoidance).  But then if $J\neq 0$ $s$ must be in an
embedded prime, contrary to our assumption.
\end{proof}

\begin{lemma}\label{eleminkl}
There are inclusions
\begin{itemize}
\item $\I_i\I_j\subset \I_{i+j+1}$ and
\item $\I_X\I_j \subset \I_{j+1}$.
\end{itemize}
\end{lemma}
\begin{proof}
Note that the second inclusion is the special case $i=0$ of the
first one.  Let $C_i \subset X$ be the union of the embedded components of
$Y \bigcap X^{(i)}$, and let $B$ be the union of the $C_i$ and the
$D_i$ (from Proposition \ref{cmcod1}).  For all $x \notin B$ we have equations
\[ \I_{i,x}=\I_{Y,x}+\I_{X,x}^{i+1},\]
\[ \I_{j,x}=\I_{Y,x}+\I_{X,x}^{j+1},\]
\[ \I_{i+j+1,x}=\I_{Y,x}+\I_{X,x}^{i+j+2}.\]
By multiplying the two first we get
\begin{eqnarray*}
\I_{i,x}\I_{j,x}=\I_{Y,x}\I_{Y,x}+\I_{Y,x}\I_{X,x}^{j+1}+\I_{Y,x}
\I_{X,x}^{i+1}+\\
+\I_{X,x}^{i+1}\I_{X,x}^{j+1} \subset
\I_{Y,x}+\I_{X,x}^{i+j+2} =\I_{i+j+1,x}.
\end{eqnarray*}

The quotient 
\[(\I_i\I_j+\I_{i+j+1}/\I_{i+j+1})\]
is thus an Ideal on $Y_{i+j+1}$ supported on $B$, and is therefore
zero by the previous lemma.
\end{proof}

Because of the inclusions $\I_X\I_j \subset \I_{j+1} \subset \I_j$ we
have another inclusion
\[\I_{j+1}/\I_X\I_j \subset \I_j/\I_X\I_j,\]
which is an inclusion of $\O_X$-Modules.  Define the $\O_X$-Module
$\L_j$ as the cokernel of this map, so that we have a short exact sequence
\begin{equation}
\label{defavL_j}
0 \rightarrow \I_{j+1}/\I_X\I_j \rightarrow \I_j/\I_X\I_j
\rightarrow \L_j \rightarrow 0.
\end{equation}
\begin{prop}
\label{strukknippl}
 There is an exact sequence of $\O_Z$-Modules
\[0\rightarrow i_{\ast}(\L_j) \rightarrow \O_{Y_{j+1}} \rightarrow
\O_{Y_j} \rightarrow 0,\]
where $i:X\hookrightarrow Z$ is the inclusion.

\end{prop}
\begin{proof}
From the short exact sequence \ref{defavL_j} we get
\[0\rightarrow \I_{j+1}\rightarrow \I_j \rightarrow i_{\ast}(\L_j)
\rightarrow 0.\]
We place this in a larger diagram of $\O_Z$-Modules
\[\xymatrix{
         &          &    & 0\ar[d]              & \\
         & 0\ar[d]  &    & i_{\ast}(\L_j)\ar[d] & \\
0\ar[r]  &\I_{j+1}\ar[r]\ar[d] & \O_Z
         \ar[r]\ar[d]^{=}&\O_{Y_{j+1}}\ar[r]\ar[d] &0\\
0\ar[r]  &\I_j\ar[r]\ar[d] & \O_Z
         \ar[r]&\O_{Y_j}\ar[r]\ar[d] &0\\
&i_{\ast}(\L_j)\ar[d]&&0&\\
&0&&&}\]
Here the Snake Lemma shows that we get $i_{\ast}(\L_j)$ in the upper
         right corner.
\end{proof}

\begin{example}[Removing embedded components]
Consider a linear subspace $\PP^n\subset \PP^{n+2}=\Proj
k[z_0,...,z_n,x,y]$ with ideal $(x,y)$.  For a linear form
$L=a_0z_0+...+a_nz_n$, there is 
a multiplicity four structure with ideal $(x^2+Ly,y^2)$.  Intersecting
this with the second infinitesimal neighbourhood of $\PP^n$ we
get $(x^2+Ly,y^2,x^3)$.  This ideal has an embedded component, a
$\PP^{n-1}$ with ideal $(L,x,y)$.  Removing this embedded component
gives us the ideal $(x^2+Ly,xy,y^2)$.  Intersecting this with the
first infinitesimal neighbourhood gives the ideal $(Ly,x^2,xy,y^2)$,
which again has an embedded component.  Removing this gives us
$(x^2,y)$.  Thus the ideals of the terms of the $S_1$-filtration are
\[(x,y)\supset (x^2,y) \supset (x^2+Ly,xy,y^2) \supset (x^2+Ly,y^2).\]

\end{example}

\begin{bem}
There are known criteria for when a double structure arises from the
Ferrand construction.  See for instance Manaresi,\cite{Man84}, where one can
find the following example:  If the original scheme $X$ is given, in
projective space, by the ideal $(x^6,y^6)$, the scheme with ideal
$(x^9,y^8)$ has twice the multiplicity of $X$ in each point but its ideal does
not contain the monomial $x^6y^6$, and thus it does not contain the ideal
$(x^6,y^6)^2$.  This shows that there are, for some schemes,
multiple structures that cannot be obtained via the Ferrand
construction, even in the case of complete intersections.  Also,
reduced schemes with bad singularities can render the Ferrand construction
insufficient.  Our smoothness assumption is included to remove these
problems, though much can be said in a more general setting.
\end{bem}

If the ambient space is $Z=\PP^N$, or if a projective embedding
of $Z$ is fixed, one can compute Hilbert polynomials:

\begin{prop}
\label{hilbofmult}
Assume $Y$ satisfies $(S_1)$.  Then the Hilbert polynomial of $Y$ can
be computed from the $S_1$-filtration as follows:

\[\Hilb(Y)=\Hilb(X)+\sum_j \Hilb(\L_j).\]
\end{prop}

\begin{proof}
This follows readily from the exact sequences

\[0 \rightarrow i_{\ast}(\L_j)\rightarrow \O_{Y_{j+1}} \rightarrow
\O_{Y_j} \rightarrow 0,\]
the fact that
\[\Hilb \L_j=\Hilb i_{\ast}\L_j,\]
that $Y=Y_k$ (see Remark \ref{YS1})
and the additivity of the Hilbert polynomial.
\end{proof}

\begin{bem}
The basic theory presented in this section was first considered by
B\v{a}nic\v{a} and Forster in \cite{BF81} for curves, later generalized by
Holme in \cite{Holpre} to higher dimensions.  Manolache has also
considered similar constructions.  The $S_1$-filtration was
considered by B\v{a}nic\v{a} and Forster for curves (then it is a
Cohen-Macaulay filtration).  They proved all results applicable to
curves (except the computation of Hilbert polynomials) in that case.
The idea to consider an $S_1$-filtration in higher dimensions is due
to Holme.  Also, the statements and proofs
of Lemmas \ref{alglemma} and \ref{eleminkl} are taken from Holme's
preprint \cite{Holpre}.  Compare also with the similar filtrations
considered by Manolache, e.g. in \cite{Man94}.
\end{bem}

\subsection{Completions}
In this section we will prove some results linking Cohen-Macaulay
schemes and locally free sheaves.  Because Cohen-Macaulay is a local
condition, we can consider our schemes point by point, and use the
properties of completion.\\

Let $X\subset Y\subset Z$ be as before.  From Sequence \ref{defavL_j}
and Proposition \ref{strukknippl} there are short exact sequences
linking the $\O_X$-Modules $\L_j$ to the structure sheaves and Ideals
of terms of the $S_1$-filtration.  We do not have an
$\O_X$-Module structure on the terms $\O_{Y_j}$, even though they are all
supported on $X$.  On the other hand we can do the following:\\

\begin{itemize}
\item Localize to any point $x\in X$.
\item Complete with respect to the maximal ideal $\mathfrak{m}_x\subset
  \O_{Z,x}$. 
\end{itemize}

This gives an ``infinitesimal Module structure'' to these structure
sheaves.\\

\begin{bem}
In the classical topology, there are Module structures on
multiple structures, coming from projections of tubular
neighbourhoods.  We are in a sense mimicking the results one might get
from that construction.  These techniques were used by B\v{a}nic\v{a}
and Forster in \cite{BF81} for curves.
\end{bem}

 Since $X$ and $Z$ are smooth, the structure map from $\widehat{\O_{Z,x}}$
 to $\widehat{\O_{X,x}}$ admits a retraction, thus giving any
 $\widehat{\O_{Z,x}}$-module a structure of
 $\widehat{\O_{X,x}}$-module.
 By choosing a suitable set of coordinates, we can assume that

\[ \widehat{\O_{Z,x}}\iso k[[X_1,...,X_N]],\]

\[ \widehat{\O_{X,x}}\iso k[[X_1,...,X_n]].\]
  For any scheme $Y$ such that $X\subset Y\subset X^{(k)}$, the
  algebra $\widehat{\O_{Y,x}}$ becomes a finite
  $\widehat{\O_{X,x}}$-module, being a quotient of
  $k[[X_1,\dots,X_N]]/I^{k+1}$, where $I$ is the ideal $(X_{n+1},\dots,X_N)$.

By completion we have short exact sequences

\[ 0 \rightarrow \widehat{(\I_{j+1}/\I_j\I)_x} \rightarrow
  \widehat{(\I_j/\I_j\I)_x}
\rightarrow \widehat{(\L_j)_x} \rightarrow 0\]

\noindent and

\[0\rightarrow \widehat{i_{\ast}(\L_j)_x}\rightarrow \widehat{(\O_{Y_{j+1}})_x}
\rightarrow \widehat{(\O_{Y_j})_x}\rightarrow 0\] 
\noindent of $\widehat{\O_{X,x}}$-modules, the latter consisting of
$\widehat{\O_{X,x}}$-modules by virtue of the above lemma.

\begin{lemma}
\label{CM<-->fri}
Let $Y\subset Z$ be a scheme with $Y_{\red}=X$.  Then $Y$ is
Cohen-Macaulay if and only if $\widehat{\O_{Y,x}}$ is a free
$\widehat{\O_{X,x}}$-module for all points $x\in X$, of constant rank.
\end{lemma}
\begin{proof}
The first remark is that a local ring is Cohen-Macaulay if and only if its
completion is.  Now $\widehat{\O_{X,x}}$ is a regular local
ring, and the algebra $\widehat{\O_{Y,x}}$ is a finite
$\widehat{\O_{X,x}}$-module.  In this situation we have the following
equivalences: $\widehat{\O_{Y,x}}$ is Cohen-Macaulay $\Longleftrightarrow$
$\widehat{\O_{Y,x}}$ is a free $\widehat{\O_{X,x}}$-module $\Longleftrightarrow$
$\widehat{\O_{Y,x}}$ is a projective $\widehat{\O_{X,x}}$-module
$\Longleftrightarrow$ the 
multiplicity of $\widehat{\O_{Y,x}}$ is constant over all the prime ideals
of $\widehat{\O_{X,x}}$.  The first equivalence is \cite{Eis95} Corollary
18.17 (which applies since $\widehat{\O_{X,x}}$ is considered as a
subring of $\widehat{\O_{Y,x}}$), the second is because $\widehat{\O_{X,x}}$ is a regular local ring, the third is
\cite{Eis95} Exercise 20.13.  The part concerning constant rank is
deduced because we can localize in the zero ideal of
$\widehat{\O_{X,x}}$, complete it,
and the resulting ring is independent of the chosen point $x$, being
the completion of the local ring of the generic point.
\end{proof}

\begin{bem} This constant rank is equal to the multiplicity of $Y$.
  Indeed, we can define the multiplicity of any nilpotent scheme $Y$
  with support $X$ to be the multiplicity of the Artin ring
  $\O_{Y,\xi}$, where $\xi$ is the generic point of $X$.  This number
  is then equal to the rank of $\widehat{\O_{Y,\xi}}$ as an
  $\widehat{\O_{X,\xi}}$-module.
\end{bem}

\begin{prop}
\label{CMlokfri}
\begin{itemize}

\item[a)] If $Y_j$ is Cohen-Macaulay and $\L_j$ is locally free on $X$,
  then $Y_{j+1}$ is Cohen-Macaulay.
\item[b)] If $Y_j$ and $Y_{j+1}$ are Cohen-Macaulay, then $\L_j$ is
  locally free on $X$.
\end{itemize}
\end{prop}
\begin{proof}
For $a)$, let $x\in X$ be any point.  We need to show that
$\O_{Y_{j+1},x}$ is a Cohen-Macaulay ring.  A local ring is
Cohen-Macaulay if and only if its completion is, so we are reduced to
proving that $\widehat{\O_{Y_{j+1},x}}$ is Cohen-Macaulay.  Now 
\[0\rightarrow \widehat{i_{\ast}(\L_j)_x}\rightarrow \widehat{(\O_{Y{j-1}})_x}
\rightarrow \widehat{(\O_{Y_j})_x} \rightarrow 0\]
is a sequence of $\widehat{\O_{X,x}}$-modules.  The leftmost is assumed
free, the last one is free by Lemma \ref{CM<-->fri}.  Thus the
middle one is free also, and we can invoke the lemma again to
conclude (note that the rank of all modules involved can be measured
by localizing in the zero ideal.  Thus we see that it is independent
of the choice of the point $x$).\\

For part $b)$ the same exact sequence and the same lemma show that the
rank of the sheaf $\L$ is independent of the point $x\in X$,
and that it is a free module in each point of $X$.  Thus it is locally
free. 
\end{proof}

\begin{bem} This central proposition can also be proved without
  completions, using homological properties of Cohen-Macaulay rings.
  For details, see \cite{Vat01}.
\end{bem}

\begin{prop}
\label{s1torfri}
The $\widehat{\O_{X,x}}$-modules $\widehat{\O_{Y_j,x}}$ are torsion
free.
\end{prop}
\begin{proof}
Let $A=\widehat{\O_{X,x}},B=\widehat{\O_{Z,x}},
C=\widehat{\O_{Y_j,x}}, I_j=\widehat{(\I_{j})_x},I=\widehat{(\I_X)_x}$.
  Suppose there is a relation $ac=0, 0\neq a\in A, c\in C$.  $A$ can
  be considered as a subring of $B$, and then $a\notin I$.  We need to
  show that $c=0$.  Choose any element $b\in B$ that represents $c$, i.e.
  $c=b \mod I_j$.  Then clearly
\[ac =0 \Longleftrightarrow ab \in I_i\]
and
\[ab\in I_i, a\notin I \Longrightarrow b\in I_i\]
since $I_i$ is $I$-primary.  Thus $c=b\mod I_i =0$.

\end{proof}

\begin{cor}
The $\O_X$-Modules $\L_j$ are torsion free.
\end{cor}
\begin{proof}
$\widehat{((\L_j)_x)}$ is a submodule of $\widehat{\O_{Y_j,x}}$, and is
therefore torsion free.  This means that the canonical map
$\widehat{((\L_j)_x)} \rightarrow \widehat{((\L_j)_x)^{\vee \vee}}$ is
injective.  Since completion is exact and faithful and commutes with
taking homomorphism groups the map
$(\L_j)_x\rightarrow (\L_j)_x^{\vee \vee}$ is injective too.
Thus $(\L_j)_x$ is torsion free.
\end{proof}

We end this section with a theorem that we will use a lot in
applications of this theory.  We give a partial order on subschemes of
$Z$ by inclusion, and a partial order on quotients of a sheaf $\F$ by the
following rule: $\F\surj \K_1\surj\K_2 \Leftrightarrow \K_2\leq\K_1$.
\begin{thm}
\label{ekvlokfrilcm}
Let $Y\subset Z$ be a Cohen-Macaulay scheme with reduced subscheme
$X$.  There is a $1-1$ correspondence between locally free quotient
$\O_X$-Modules of $\I_Y/\I_Y\I_X$ and Cohen-Macaulay schemes $W$
satisfying $\I_Y\I_X\subset \I_W \subset \I_Y$.  This correspondence
preserves the partial orders defined above.  If there are surjections
$\I_Y/\I_Y\I_X\surj \K_1 \surj \K_2$ and $W_i$ corresponds to $\K_i$
then there is a surjection $\I_{W_2}/\I_{W_2}\I_X\surj (\ker(\K_1\surj
\K_2))$.
\end{thm}
\begin{proof}
The first part can be proven along the same lines as Proposition
\ref{CMlokfri}, since the fact that the schemes involved are part of an
$S_1$-filtration is not essential.\\

The second part follows from the first part if we show
$\I_Y/\I_Y\I_X\surj \K_1\surj\K_2\Rightarrow W_2\subset W_1$, where
$W_i$ is the scheme corresponding to $\I_Y/\I_Y\I_X\surj \K_i$.  This
can easily be seen from the diagram of $\O_X$-Modules
\[\xymatrix{
&&&0\ar[d]&\\
&0\ar[d]&&\M\ar[d]&\\
0\ar[r]&\I_{W_1}/\I_Y\I_X\ar[r]\ar[d]&\I_Y/\I_Y\I_X\ar[r]\ar@{=}[d]&\K_1\ar[r]\ar[d]
&0\\
0\ar[r]&\I_{W_2}/\I_Y\I_X\ar[r]\ar[d]&\I_Y/\I_Y\I_X\ar[r]&\K_2\ar[r]\ar[d]&0\\
&\M\ar[d]&&0&\\
&0&&&}\]
Here $\M$ can be defined as the kernel of the right hand map, and the
Snake Lemma shows that $\M$ also appears in the lower left corner.
Thus the partial orders are preserved.\\

The last part can be deduced from the diagram above by considering the
sequence of isomorphisms
\[\M\iso \frac{\I_{W_2}/\I_Y\I_X}{\I_{W_1}/\I_Y\I_X}\iso
\frac{\I_{W_2}}{\I_{W_1}}\iso
\frac{\I_{W_2}/\I_{W_2}\I_X}{\I_{W_1}/\I_{W_2}\I_X}.\] 

\noindent The first isomorphism is the Snake Lemma applied to the diagram
above, the second is the classical Second Isomorphism Theorem.  So is
the third, noting that $\I_{W_2}\I_X\subset \I_{W_1}$, which is the kernel
of the surjection $\I_{W_2}\surj\M$.
\end{proof}

\begin{bem}
\label{spesforproj}
Note that many of the properties proven in this section can be
deduced more easily in case $X$ is a linear subspace of $Z=\PP^N$.  In
that case the projection from $Z$ minus a complementary subspace onto
$X$ restricts to a map from each (embedded) multiple structure on $X$
onto $X$, and the various structure sheaves and quotients of such
become $\O_X$-Modules.  This can of course happen for other types of
multiple structures as well.  Thus it is natural to make the following
definition.
\end{bem}

\begin{defn}
Generalizing \cite{BE95} (page 723), we say that a multiple structure $Y$
on $X$ is a {\em split } structure if the canonical inclusion
$X\hookrightarrow Y$ admits a retraction $Y\rightarrow X$.  By
restriction, all terms $Y_k$ of the $S_1$-filtration will then also
satisfy this condition. 
\end{defn}

\begin{example}
Take any variety $X$ and any vector bundle $\E$, form the total space
$E$ (or its projective completion)
and embed $X$ as the zero section.  Any multiple structures on
$X$ in this embedding is split, the retraction being the restriction
of the projection map $E\rightarrow X$.  For instance, all
infinitesimal neighbourhoods of $X$ in this embedding split.
\end{example}

We end this chapter by computing some dualizing sheaves. Because of
the smoothness conditions, the ambient space $Z$ has the dualizing
sheaf $\omega_Z$, which is equal to the canonical sheaf, and similarly
for $X$.  Recall that (with $c:=\codim (X,Z)$)
\[\omega_X\iso {\mathcal Ext}_{\O_Z}^c(\O_X,\omega_Z),\]

\noindent and if $Y$ is a Cohen-Macaulay scheme whose reduced
subscheme is $X$, then
\[\omega_X\iso{\mathcal Hom}_{\O_Y}(\O_X,\omega_Y)\]
(see for instance \cite{AK70}).  Also, since $X$ is a locally complete
intersection, we have an isomorphism
\[\omega_X\iso \wedge^c \N_{X/Z}\otimes_{\O_Z}\omega_Z\]
and similarly for $Y$ in the case that $Y$ is a locally complete
intersection.  The following result on canonical Modules ties these
together in the simplest case.  This result was also presented by
B\v{a}nic\v{a} and Forster for curves, and by Manolache in general.  

\begin{prop}
\label{dualiknippe}
Suppose $Y$ is a double structure on $X$, and a locally complete
intersection.  Let the thickening $X\subset Y$ be given by the
sequence
\[0\rightarrow \I_Y/\I_X^2\rightarrow \I_X/\I_X^2\rightarrow
\L\rightarrow 0\]
for an invertible sheaf $\L$.  Then there is an isomorpishm
\[\omega_Y|_X\iso \omega_X\otimes \L^{-1}.\]
\end{prop}

\begin{proof}
Note that
\[\I_Y/\I_Y^2|_X\iso \I_Y/\I_Y\I_X\]
and consider the sequence

\[0\rightarrow \I_X^2/\I_Y\I_X\rightarrow\I_Y/\I_Y\I_X\rightarrow
\I_Y/\I_X^2 \rightarrow 0.\] 
Splicing this to the sequence mentioned in the statement, and noting
that $\L^2\iso (\I_X/\I_Y)^2\iso \I_X^2/\I_X\I_Y$ we get an exact
sequence
\[0\rightarrow \L^2\rightarrow \I_Y/\I_Y\I_X\rightarrow
\I_X/\I_X^2\rightarrow \L\rightarrow 0.\] Taking the top exterior
power (or {\em determinant}) of this sequence, we get
\[\L^2\otimes \wedge^c\I_X/\I_X^2\iso \wedge^c \I_Y/\I_Y\I_X\otimes
\L\] which can be rewritten
\[\L\otimes \wedge^c\N_{Y/Z}|_X\iso \wedge^c\N_{X/Z}.\]
Combining this with the two isomorphisms
\[\omega_X\iso \wedge^c \N_{X/Z}\otimes_{\O_Z}\omega_Z\]
and
\[\omega_Y\iso \wedge^c \N_{Y/Z}\otimes_{\O_Z}\omega_Z\]
we get
\[\omega_Y|_X\iso \omega_X\otimes \L^{-1}.\]

\end{proof}

\begin{bem}
B\v{a}nic\v{a} and Forster, as well as Holme, considered the notion of
{\em primitive} structure.  Basically, this is the following:  assume
$X$ is given (locally) by $m=\codim X$ equations $f_1,\dots,f_m$.
Then a primitive multiple structure $Y$ on $X$ is given (locally) by
the equations $f_1^l,f_2,\dots,f_m$.  If the double substructure
given by $f_1^2,f_2,\dots,f_m$ corresponds to a sequence

\[0\rightarrow \J/\I_X^2\rightarrow \I_X/\I_X^2\rightarrow
\L\rightarrow 0\]
then
\[\omega_Y|_X\iso \omega_X\otimes \L^{-(l-1)}.\]

Also, each $\L_j\iso \L^j$.  See \cite{BF81},\cite{BF86} and \cite{Holpre}.
\end{bem}


\section{\label{codimtwo} Classification of multiple structures of low degree on linear subspaces of codimension two}

In this chapter we present a list of embedded multiple structures
supported on a linear subspace $\PP^n\subset \PP^{n+2}$.  Since the
case of curves is 
very special, we will throughout assume $n\geq 2$.  The list will
contain all type I multiple structures up to and including multiplicity
(=degree) five.  For structures not of type I, the classification will
be completed in the final section, up to degree four.  Up to multiplicity four
all structures were found by Manolache \cite{Man92}, and (with a small
mistake) by Holme \cite{Holpre}.  The lists in this section where
found in my Hovedfagsoppgave (Master's thesis), \cite{Vat98}, except for
some small mistakes.\\  

We will fix the following notation in this chapter.

\begin{notation}
The fixed ambient space is $\PP^{n+2}=\Proj k[z_0,\dots,z_n,x,y]$.  We
denote by  $X:=\PP^n=\Proj k[z_0,\dots,z_n]\subset \PP^{n+2}$ the linear
subspace that is the support of the multiple structures, with ideal
$I=(x,y)\subset S:=k[z_0,\dots,z_n,x,y]$ and Ideal $\I \subset
\O_{\PP^{n+2}}$.  
\end{notation}

Now I wish to recall the exact statement we need for our inductive
argument.  It can be extracted from Proposition \ref{CMlokfri} and
Theorem \ref{ekvlokfrilcm}.

\begin{prop}
\label{induksjon}
Let $Y$ be a multiple structure, $Y_{\red}=X$.  Assume that all the $Y_k$
are Cohen-Macaulay.  Then all the $\L_k$ are locally free sheaves, and we
have short exact sequences
\[0\rightarrow \J_{k+1}/\I\J_k\rightarrow
\J_k/\I\J_k\rightarrow \L_k\rightarrow 0\]
for all $k$.
\end{prop}

Thus, in order to construct multiple structures satisfying the
conditions of the proposition, we can proceed inductively from
multiple structures of lower multiplicity.  For each $k$, all we need
for the induction is to know which locally free quotients
$\J_k/\I\J_k$ have.\\

\begin{defn}
Recall that Cohen-Macaulay schemes satisfying the conditions of the
proposition are said to be of {\em type I}.  Throughout this chapter
the term ``multiple structures'' refers to Cohen-Macaulay multiple
structures of type I.
\end{defn}

\begin{bem}A note about the list: we present representatives for the various
projective equivalence classes, i.e. up to coordinate changes.  Some
of the representatives include an ``undetermined'' polynomial of a given
degree (see e.g. the polynomial $G$ occuring in Proposition \ref{primitiv}).  By abuse of the term
``projective equivalence'' we don't differ between various classes of
these polynomials.  Furthermore,
we will use a modified version of the above proposition.  Since a
thickening corresponding to an $\L$ which is an extension of locally
free sheaves $\E,\F$ of lower rank is a composition of thickenings
corresponding to the sheaves $\E,\F$, we will only consider $\L$s that
are not extensions.  In our computations this usually ensures that
each $\L$ is an invertible sheaf.

\end{bem}

\subsection{Multiplicity one, two and three}

In multiplicity one, there is only one multiple structure $Y$ with
$Y_{\red}=X$, namely $Y=X$.\\

Using Proposition \ref{induksjon}, we want exact sequences

\[0\rightarrow \J/\I^2\rightarrow \I/\I^2\rightarrow \L\rightarrow 0\]
\noindent
where $\L$ is an invertible sheaf.  Note that $\I/\I^2\iso
\O_X(-1)^{\oplus 2}$ with ``generators'' $x$ and $y$.  Since $n\geq
2$ $\L$ has to be isomorphic to $\O_X(-1)$.  Thus $\J/\I^2\iso
\O_X(-1)$, and $J$ is generated modulo $I^2$ by a polynomial of degree
one.  Up to linear equivalence, this can be taken to be $y$.\\

There is thus, up to linear equivalence, only one double structure $Y$ with
$Y_{\red}=X$, and its ideal can be taken to be $(y)+I^2=(x^2,y)$.\\

For multiplicity three, we will use a general statement.  Letting
$\nu=2$ in the following proposition gives the multiplicity three
case.

\begin{prop}
\label{primitiv}
Let $Y$ be the $\nu-$multiple structure with ideal $K=(x^{\nu},y)$, where
$\nu\geq 2$.  There are, up to linear equivalence, three structures of
multiplicity $\nu +1$ containing $Y$.  Their ideals can be taken to be
\begin{itemize}
\item $J_1=(x^{\nu+1},y)$
\item $J_2=(x^{\nu},xy,y^2)$
\item $J_3=(x^{\nu}+Gy,xy,y^2)$
\end{itemize}
where $G\in k[z_0,\dots,z_n]_{\nu-1}$.
\end{prop}

\begin{proof}
Note that $\K/\I\K\iso \O_X(-\nu)\oplus \O_X(-1)$ with ``generators''
$x^{\nu}$ and $y$.  Thus there are only two possibilities for the
invertible sheaf $\L$ in the short exact sequence 
\[0\rightarrow \J/\I\K\rightarrow \K/\I\K\rightarrow \L\rightarrow 0,\]
namely $\L\iso \O_X(-\nu)$ and $\L\iso \O_X(-1)$.  In the first case,
$\J/\I\K\iso \O_X(-1)$ and thus $J$ is generated modulo $IK$ by a
linear polynomial contained in $K$.  This way we get the ideal $J_1$.
In the second case $\J/\I\K\iso \O_X(-\nu)$, and $J$ is generated
modulo $IK$ by a polynomial of degree $\nu$, contained in $K$.  Call
this polynomial $F$.  Then $F=ax^{\nu}+Gy$, where $a$ is a scalar and
$G$ is a polynomial of degree $\nu-1$.  Note that since we include all
of $IK$ in the ideal $J$ we can assume that $G\in k[z_0,\dots,z_n]$.
If $a\neq 0$ we can rescale $x$ to get $F=x^{\nu}+Gy$, and we get the
two ideals $J_2,J_3$ depending on whether $G$ is zero or not.\\

Consider the case $a=0$, so $F=Gy$.  Let $G=\Pi G_i^{g_i}$ be the
factorization of $G$ in irreducible polynomials.  For each $i$ we have 

\[(J:yG/G_i^{g_i})=(x,y,G_i^{g_i})\]
\noindent
which is a primary ideal for the prime ideal $(x,y,G_i)$.  This is
thus an embedded associated prime ideal for $(F)+IK$.  Thus this ideal
is not Cohen-Macaulay, and therefore excluded from our list.
\end{proof}

\begin{bem}

There are a couple of remarks that should be made concerning the
case $a=0$ in the above proof.  One is that the problem can be spotted
directly from the sequence

\[0\rightarrow
\O_X(-\nu)\stackrel{\left(\begin{array}{c}0\\G\end{array}\right)}{\rightarrow}
\O_X(-\nu)\oplus \O_X(-1) \rightarrow \L\rightarrow 0.\] 
Since the left map is the zero map from $\O_X(-\nu)$ to $\O_X(-\nu)$
and is multiplication by $G$ from $\O_X(-\nu)$ to $\O_X(-1)$, the only
way such a sequence can be exact is if $\L$ is a copy of $\O_X(-\nu)$
and in addition a torsion part, namely the structure sheaf of the zero
scheme of $G$, twisted by $-1$.  In particular, the left hand map does
not occur in such a short exact sequence with $\L$ locally free.\\

The other remark to be made is that one might consider removing the embedded components.  Then the generator $Gy$ in the ideal has to be replaced by
$y$ alone, and we find again the ideal $J_1$.  This is no surprise;
from the first remark we know that in this case $\L$ would contain a
copy of $\O_X(-\nu)$.  This piece alone will give the Ideal $\J_1$.
In addition, $\L$ has a torsion part which will add some embedded
components to the ideal $J_1$.
\end{bem}

Summarizing, we get

\begin{thm}
The type I structures on $X$ of multiplicity less than or equal to
three are, up to linear equivalence:

\begin{eqnarray*}
(x,y)\\
(x^2,y)\\
(x^3,y)\\
(x^2,xy,y^2)\\
(x^2+z_0y,xy,y^2)
\end{eqnarray*}
(the linear polynomial $G$ can be taken to be $z_0$ in the last
structure).
\end{thm}

\subsection{Multiplicity four and five}
Most of what follows can be constructed using exactly the same
techniques as we used in the last section.  There are some
additional complications arising, and these are the only ones we will
give proofs of.  Full proofs can be found in \cite{Vat98}.  A rather
interesting non-existence result is postponed to the next section.\\

\begin{thm}
The type I structures on $X$ of multiplicity four are, up to linear
equivalence, the following:
\begin{eqnarray*}
(x^4,y)\\
(x^2,y^2)\\
(xy,x^2+y^2)\\
(x^3,xy,y^2)\\
(F_2y^2-F_3xy,F_1y^2-F_3x^2,F_1xy-F_2x^2,x^3,x^2y,xy^2,y^3)\\
(x^3+Gy,xy,y^2)\\
(x^2+z_0y,y^2)\\
(x^2+xy+z_0y,y^2)\\
\end{eqnarray*}
\noindent
Here $(F_2y^2-F_3xy,F_1y^2-F_3x^2,F_1xy-F_2x^2,x^3,x^2y,xy^2,y^3)$
only occurs if $\dim X=2$, in which case the polynomials $F_i\in
k[z_0,z_1,z_2]_s$, $s\geq 1$ satisfy $\bigcap V_+(F_i)=\emptyset$.
$G$ is a polynomial of degree 2 in $k[z_0,\dots,z_n]$.   The pairs of 
structures $(x^2,y^2)$, $(xy,x^2+y^2)$ and
$(x^2+z_0y,y^2)$,$(x^2+xy+z_0y,y^2)$ are linearly equivalent if the
characteristic is different from two.
\end{thm}

\begin{bem}
Thus there are two new phenomena:  we get more structures for low
dimension (with an infinite number of different Hilbert polynomials),
and the special features of characteristic two show up in the list.
\end{bem}

First we will see how the characteristic enters into the analysis.

\begin{prop}
Starting with the multiple structure $Y$ with ideal $K=(x^2,xy,y^2)$,
there are the following possibilities, up to linear equivalence, for a
structure of multiplicity four (assuming $\dim X\geq 3$; there are
additional structures if $\dim X=2$):

\begin{eqnarray*}
(x^3,xy,y^2)\\
(x^2,y^2)\\
(x^2+y^2,xy)
\end{eqnarray*}
where the last two are linearly equivalent unless the characteristic
is two.
\end{prop}
\begin{proof}

Note first that 
\[\K/\I\K\iso \O_X(-2)^{\oplus 3}\]
with ``generators'' $x^2,xy,y^2$.  The only way to get a surjection
from this sheaf onto an invertible sheaf $\L$ is to choose $\L\iso
\O_X(-2)$ (this is where the assumption $\dim X\geq 3$ is needed).
Thus, in the short exact sequence
\[0\rightarrow \J/\I\K \rightarrow \K/\I\K\rightarrow \L \rightarrow
0\]
\noindent
the term $\J/\I\K\iso \O_X(-2)^{\oplus 2}$; thus $J$ is generated, modulo $IK$,
by two homogeneous polynomials of degree two, contained in $K$.  Let
these polynomials be
\[F_1:=a_1x^2+b_1xy+c_1y^2, \mbox{  }F_2:=a_2x^2+b_2xy+c_2y^2.\]
\noindent
By replacing $F_1$ with $F_1-b_1/b_2F_2$ (or renumbering the $F_i$ in
case $b_2=0$) we can assume that $b_1=0$ (we are only interested in
the ideal generated).\\

Now if both $a_1=a_2=0$, we get $F_1=c_1y^2$ and $F_2=b_2xy+c_2y^2$.
Normalizing to $c_1=b_2=1$ and replacing $F_2$ with $F_2-c_2F_1$ we
get the ideal $(xy,y^2)+IK=(x^3,xy,y^2)$.\\

Assume $a_1=1$.  Replacing $F_2$ by $F_2-a_2F_1$ we can assume
$a_2=0$.  There are two cases to consider, namely $c_1=0$ and $c_1\neq
0$ (in which case we can assume $c_1=1$ be rescaling $y$).  These two
cases turn out to give basically the same analysis, so we restrict our
attention to the case $c_1=0$.\\

Thus we consider $F_1=x^2$ and $F_2=b_2xy+c_2y^2$.  If either of
$b_2,c_2$ is zero we get ideals $(x^2,xy,y^3)$ and $(x^2,y^2)$.  So
assume $b_2\neq 0 \neq c_2$.  Rescaling $y$ we can assume $c_2=1$.
Then, rescaling $x$, we can assume $b_2=1$:  now $F_1=x^2,F_2=xy+y^2$.
If the characteristic is different from two, the ideal generated by
$F_1$ and $F_2$ can be put in the form $(x^2,y^2)$ through a linear change of
coordinates.  Explicitly, let
\[
\left(\begin{array}{c}z\\w\end{array}\right)=\left(\begin{array}{cc}1&0\\1/2&1
  \end{array}\right) \left(\begin{array}{c}x\\y\end{array}\right)\]
\noindent
Then $z^2=x^2$ and $w^2=1/4x^2+xy +y^2$, and the ideal generated by
these two squares is obivously the same as the ideal generated by
$F_1$ and $F_2$.  If the characteristic is two, however, this cannot
be done.  It is most easily seen if we try going the other way:
starting from $z^2,w^2$ and making a linear change of coordinates, we
get e.g. $z^2=(ax+by)^2=a^2x^2+b^2y^2$.  Thus we can in no way get the
monomial $xy$ as part of our generating polynomial.  To get $(F_1,F_2)$
in the form from the statement, use $x\mapsto x-y,y\mapsto y$.\\

The remaining choices of coordinates can easily be turned into one of
the situations discussed above by quite obvious linear changes of coordinates.
\end{proof}

The next new phenomenon that enters is the presence of torsion in the
sheaves $\K/\I\K$.  The following proposition is a continuation of a
part of Proposition \ref{primitiv}:

\begin{prop}
Let $Y$ be the multiple structure on $X$ with ideal
$K=(x^{\nu}+Gy,xy,y^2)$, where $\nu\geq 2$.  Then, in order to have a
surjection from
$\K/\I\K$ onto a locally free sheaf $\L$, it is necessary and sufficient that
$\L$ is isomorphic to a direct summand of $\O_X(-\nu)\oplus
\O_X(-2)\iso \coker (\O_X(-2)\stackrel{y^2}{\rightarrow} \K/\I\K)$.  In
particular, $y^2$ is always contained in the ideal of a simple
thickening of $Y$.
\end{prop}
\begin{proof}
Note first that 
\[IK=(x^{\nu +1}+Gxy,x^{\nu}y+Gy^2,x^2y,xy^2,y^3)=(x^{\nu
  +1}+Gxy,Gy^2,x^2y,xy^2,y^3).\] 
\noindent
The term $Gy^2$ tells us that the quotient sheaf $\K/\I\K$ has
torsion.  It is easy to write down a free resolution of this sheaf:
\[\xymatrix{
0\ar[r]&\O_X(-\nu-1)\ar[r]^{\left(\begin{array}{c} 0 \\ 0 \\ G \end{array}\right)}
&\O_X(-\nu)\oplus \O_X(-2)^{\oplus 2} \ar[r] &\K/\I\K \ar[r] & 0}\]
\noindent
where the right hand map is the generator map.\\

Now giving a map from $\K/\I\K$ is equivalent to giving a map from
$\O_X(-\nu)\oplus \O_X(-2)^{\oplus 2}$ that composed with
$\O_X(-\nu-1)\rightarrow \O_X(-\nu)\oplus \O_X(-2)^{\oplus 2}$ is
zero.  If the map goes to something locally free, the entire torsion
subsheaf of $\K/\I\K$ has to be killed.  This is equivalent to
demanding that the last $\O_X(-2)$-part is in the kernel of the map,
and the map factors through the projection onto the first two summands
$\O_X(-\nu) \oplus \O_X(-2)$.  The only locally free quotients of this
sheaf are its direct summands.
\end{proof}

\begin{bem}
Note that, in this situation, $\O_X(-\nu)\oplus \O_X(-2)$ is
isomorphic to the quotient of $\K/\I\K$ by its torsion subsheaf.
\end{bem}

To give an example of how to get special structures in dimension two,
we will give a more general result which, specialized to $n=2$, gives
the example in question.

\begin{prop}
Consider $X=\PP^n\subset \PP^{n+2}$ and $Y=X^{(n-1)}$.  Then there are
infinitely many different Hilbert polynomials of thickenings of $Y$ to
a structure of multiplicity one more than the multiplicity of $Y$.
With one exception, none of these classes extends to Cohen-Macaulay
multiple structures on $\PP^{n+1}\subset
\PP^{n+3}$. 
\end{prop}
\begin{proof}
Letting $K=(x,y)^n$ be the ideal of $Y$, note that 
\[\K/\I\K\iso \O_X(-n)^{\oplus n+1}\]
with ``generators'' $x^n,x^{n-1}y,\dots,y^n$.  Now, choosing $n+1$
polynomials $F_i$ satisfying $\bigcap V_+(F_i)=\emptyset$ of the same
degree $d$, we get a short exact sequence
\[0\rightarrow \J/\I\K \rightarrow \K/\I\K\rightarrow \O_X(-n+d)\rightarrow 0\]
defining $\J$.  The $F_j$ form a regular sequence, so the image of
the map 
\[(\wedge^2(\O_X(-d)^{\oplus n+1}))(-n+d)\rightarrow
(\O_X(-d)^{\oplus n+1})(-n+d)\iso\K/\I\K\]

\noindent coming from the (twisted) Koszul
complex of the $F_i$ equals the sub-Module $\J/\I\K$.  Thus we can
easily write down generators for the ideal $J$:
\[J=(F_jx^{n-i}y^i-F_ix^{n-j}y^j)_{0\leq i< j \leq n} +I^2.\]
\noindent
Note that the Hilbert polynomial of $V(\J)$ only depends on the
degree $d$.  In fact, the Hilbert polynomial of $Y=X^{(n-1)}$ is given
by
\[\Hilb(Y)(t)=b_0+2b_1+3b_2+\dots +nb_{n-1}\]
\noindent so

\[\Hilb(V(\J))(t)=b_0+2b_1+3b_2+\dots +nb_{n-1}+b_{n-d}\]
\noindent
where 

\[b_i(t):=\binom{n+t-i}{n}=\chi(\O_X(t-i))=\Hilb(\O_X(-i))(t).\]
\noindent
In addition to the examples above, we can surject $\K/\I\K$
onto a direct summand $\O_X(-n)$.  The ideals $J$ arising in this way
can of course be extended to higher dimensional spaces.  All the
others, if extended in the obvious way, will fail to be Cohen-Macaulay
in the common zero locus of the $F_i$.
\end{proof}

\begin{thm}
The type I structures on $X$ of multiplicity five are, up to linear
equivalence, the following:
\begin{eqnarray*}
(x^5,y)\\
(x^4,xy,y^2)\\
(x^3+z_0xy,x^2y,y^2)\\
(x^3,xy,y^3)\\
(x^3,x^2y,y^2)\\
(x^3+z_0y^2,xy,y^3)\\
(F_2y^2-F_3xy,F_1y^2-F_3x_3,F_1xy-F_2x^3,x^4,xy^2,x^3y,y^3)\\
(P_1x^2+P_2xy+P_3y^2,x^3,x^2y,xy^2,y^3)\\
(x^3,x^2+y^2,y^3)\\
(x^4+G_1y,xy,y^2)\\
(x^3+G_2y,x^2y,y^2)\\
(x^3+z_0xy+G_3y,x^2y,y^2)\\
(x^2+z_0y,xy^2,y^3)\\
(x^2+z_0y+y^2,xy^2,y^3)\\
(x^3,x^2+y^2,y^3)\\
(x^2+y^2,x^2y,xy^2)\\
(x^2+xy+y^2,x^2y,xy^2)\\
(x^2+z_0y+xy+y^2,xy^2,y^3)
\end{eqnarray*}
Here the last four are special variants in characteristic two.  The
$F_i\in k[z_0,z_1,z_2]$ satisfy $\deg F_1=1+\deg F_2=1+\deg F_3\geq
2$ and $\bigcap V_+(F_i)=\emptyset$; the structure
$(F_2y^2-F_3xy,F_1y^2-F_3x_3,F_1xy-F_2x^3,x^4,xy^2,x^3y,y^3)$ only
occurs for $\dim X=2$.  The $P_i\in k[z_0,z_1,z_2]$ are all of degree
$a-2\geq 1$, and $\bigcap V_+(P_i)=\emptyset$; the structure
$(P_1x^2+P_2xy+P_3y^2,x^3,x^2y,xy^2,y^3)$ only occurs for $\dim X=2$.
The $G_i$ are polynomials in $z_0,\dots,z_n$, $\deg G_1=3$, $\deg
G_2=\deg G_3=2$.
\end{thm}

\begin{bem}
The proof of this theorem contains no new techniques compared with the
techniques introduced for lower multiplicities, with two
exceptions.  There is a multiplicity four structure that cannot be
thickened to a multiplicity five structure.  The proof of this fact
occupies the next section.  There is a multiplicity five structure
built directly on a multiplicity three structure in dimension two.
The rest of the proof of the theorem is
just tedious (but not totally trivial) calculation, and can be found
in \cite{Vat98}.
\end{bem}

\begin{prop}\label{nystruktur}  Let $Y$ be the multiple structure on
  $X=\PP^2$ with 
  ideal $(x^2,xy,y^2)$.  Then there is a multiplicity five structure
  containing $Y$, without any intermediate Cohen-Macaulay multiplicity
  four structure.  The ideal of $Y$ can be taken to be $(P_1x^2+P_2xy+P_3y^2,x^3,x^2y, xy^2,y^3)$.
\end{prop}
\begin{proof}
For any positive integer $a$ there is a short exact sequence
\[0\rightarrow \O(-a-2)\rightarrow \O(-2)^{\oplus 3}\rightarrow
\Q\rightarrow 0\]
for a bundle $\Q$ of rank two.  The left hand map is given by three
polynomials $P_1,P_2,P_3$ of the same degree $a-2$.  The image thus
corresponds to the polynomial $P_1x^2+P_2xy+P_3y^2$.  Together with the
generators of $(x,y)^3$, this polynomial generates the ideal of $Z$.
There is no intermediate multiplicity four structure because $\Q$ is
not an extension of line bundles.
\end{proof}

\begin{bem}  This proposition fills a gap from \cite{Vat98}.
\end{bem}

\begin{bem}
So far we have only seen special features in characteristic two.  The
main point was that, using linear transformations, the cross-term
$xy$ can neither appear nor disappear.  Similar problems arise in
other characteristics, but then the multiplicity must be larger.  To
show how this happens, we consider an explicit example:
\end{bem}

\begin{example}
Let $Y$ be the multiple structure on $X$ with ideal
$K=(x^3,x^2y,xy^2,y^3)$.  Then $\K/\I\K\iso \O_X(-3)^{\oplus 4}$.
Consider a short exact sequence

\[0\rightarrow \J/\I\K \rightarrow \K/\I\K \rightarrow
\O_X(-3)^{\oplus 2} \rightarrow 0.\]
\noindent
Thus the ideal $J$ will be generated, modulo $IK$, by two polynomials
of degree three.  Assume $J=(x^3, ax^2y+bxy^2+y^3)+IK$, where the
scalars $a$ and $b$ are both non-zero.  In characteristic different
from three, this can be changed to an ideal of the form $(z^3,w^3)+IK$
through a linear change of coordinates.  This cannot be done in
characteristic three.
\end{example}

\begin{bem}
The example can easily be generalized to any given characteristic $p>0$.
\end{bem}

\subsection{A non-existence result}

In this section we show that there is no multiplicity five
structure on $X$ containing the multiplicity four structure with ideal
$J:=(F_2y^2-F_3xy,F_1y^2-F_3x^2,F_1xy-F_2x^2,x^3,x^2y,xy^2,y^3)$.  The
precise result is the following:

\begin{prop}
There is no surjection 
\[\J/\I\J\rightarrow \L \rightarrow 0\]
for any invertible sheaf $\L$ on $X$.  In particular, $V(\J)$ cannot
be thickened to a Cohen-Macaulay multiplicity five structure.
\end{prop}
\begin{proof}
Note that 
\begin{eqnarray*}
IJ=(x^4,x^3y,x^2y^2,xy^3,y^4,\\
x(F_2y^2-F_3xy),x(F_1y^2-F_3x^2),x(F_1xy-F_2x^2)\\
y(F_2y^2-F_3xy),y(F_1y^2-F_3x^2),y(F_1xy-F_2x^2)).
\end{eqnarray*}
\noindent
Let us first of all calculate the rank of $\J/\I\J$.  This is
\[\dim_{K(X)}(\J/\I\J)_{\xi}\]
where $K(X)=\O_{X,\xi}$ is the function field of $X$; the local ring
of the generic point $\xi$.  Over this field the $F_i$ are invertible,
and we can define $f_1=F_1/F_3,f_2=F_2/F_3$.  There are relations
\begin{eqnarray*}
f_2xy^2=x^2y\\
f_1xy^2=x^3\\
f_1/f_2x^2y=x^3\\
f_2y^3=xy^2\\
f_1y^3=yx^2\\
f_1/f_2xy^2=yx^2
\end{eqnarray*}
\noindent
Thus the four generators $x^3,x^2y,xy^2,y^3$ only contribute with one
towards the rank of $\J/\I\J$.  The remaining three generators,
divided by $F_3$, are $f_2y^2-xy,f_1y^2-x^2,f_1xy-f_2x^2$.  They are
related by the following $K(X)$-linear relation:
\[f_1(f_2y^2-xy)+(f_1xy-f_2x^2)=f_2(f_1y^2-x^2).\]
They give a contribution of two towards the rank.  Thus,
\[\rank \J/\I\J=3.\]
\noindent
Next, we want a free presentation of this sheaf.  The first part is of
course the generator map

\[\psi:\O_X(-3)^{\oplus 4} \oplus \O_X(-s-2)^{\oplus 3}\rightarrow
\J/\I\J,\]
i.e. $\psi$ is given by
\[\psi=(x^3,x^2y,xy^2,y^3,F_2y^2-F_3xy,F_1y^2-F_3x^2,F_1xy-F_2x^2).\]
Because of the calculation above, over $K(X)$, we expect that the
kernel of $\psi$ has rank four.  This will split in two components,
one of rank three from the term $\O_X(-3)^{\oplus 4}$, and one of rank
one from $\O_X(-s-2)^{\oplus 3}$.  Precisely, let $\K:=\ker \psi$,
$\K_1:= \K\cap \O_X(-s-2)^{\oplus 3},\K_2:=\K\cap \O_X(-3)^{\oplus 4}$.
$\K_i$ can also be described as the kernels of $\psi$ composed with the
canonical injections.  Let these maps be denoted by $\psi_i$.\\

The kernel of $\psi_1$ can be represented by
\[\xymatrix{0\ar[r]&\O_X(-2s-2)\ar[r]^{\left(\begin{array}{c}F_1\\-F_2\\F_3
      \end{array}\right)} & \O_X(-s-2)^{\oplus 3}}.\]
This gives rise to an exact sequence

\[\xymatrix{0\ar[r]&\O_X(-2s-2)\ar[r]&\O_X(-s-2)\ar[r]^<<<<<{\psi_1}&\J/\I\J}.\]
\noindent
For the kernel of $\psi_2$, we have a diagram

\[\xymatrix{
0\ar[r]&\ker \psi_2\ar[r]&\O_X(-3)^{\oplus 4}\ar[r]^{\psi_2}&\J/\I\J\\
&\O_X(-s-3)^{\oplus 6}\ar[ru]^B \ar[u]&&}\]
where the $4\times 6$ matrix $B$ is given by
\[B=\left(\begin{array}{cccccc}
0 & -F_3 & -F_2 & 0 & 0 & 0\\
-F_3 & 0 & F_1 & 0 & -F_3 & -F_2\\
F_2 & F_1 & 0 & -F_3 & 0 & F_1\\
0 & 0 & 0 & F_2 & F_1 & 0   \end{array}\right) \]
\noindent
which expresses the relations
\begin{eqnarray*}
x(F_2y^2-F_3xy)\congr 0 (\mod IJ)\\
x(F_1y^2-F_3x^2)\congr 0 (\mod IJ)\\
x(F_1xy-F_2x^2)\congr 0 (\mod IJ)\\
y(F_2y^2-F_3xy)\congr 0 (\mod IJ)\\
y(F_1y^2-F_3x^2)\congr 0 (\mod IJ)\\
y(F_1xy-F_2x^2)\congr 0 (\mod IJ)
\end{eqnarray*}
\noindent
This gives an exact sequence 
\[\xymatrix{
\O_X(-s-3)^{\oplus 6} \ar[r]^<<<<{\beta}&\O_X(-3)^{\oplus
  4}\ar[r]^<<<<<{\psi_2}& \J/\I\J}.\]
\noindent
All in all we get a free presentation
\[\xymatrix{
\O_X(-s-3)^{\oplus 6} \oplus \O_X(-2s-2)\ar[r]^<<<<{\phi} &
  \O_X(-3)^{\oplus 4} \oplus \O_X(-s-2)^{\oplus 3}\ar[r]^>>>>{\psi} &
  \J/\I\J \ar[r] & 0.}\]
\noindent
$\phi$ is given by

\[\phi = \left(\begin{array}{ccccccc}
0   & -F_3& -F_2& 0   & 0   & 0   &0\\
-F_3& 0   & F_1 & 0   & -F_3& -F_2&0\\
F_2 & F_1 & 0   & -F_3& 0   & F_1 &0\\
0   & 0   & 0   & F_2 & F_1 & 0   &0\\
0   & 0   & 0   & 0   & 0   & 0   &F_1  \\
0   & 0   & 0   & 0   & 0   & 0   &-F_2  \\
0   & 0   & 0   & 0   & 0   & 0   &F_3 \end{array}\right).\]
\noindent
Now giving a surjection $\J/\I\J\rightarrow \L$, for an invertible
sheaf $\L$, is equivalent to giving a surjection $\lambda:\O_X(-3)^{\oplus
  4}\oplus \O_X(-s-2)^{\oplus 3} \rightarrow \L$ which, composed with
$\phi$, is zero.  Let $\lambda=(l_1,\dots,l_7)$.  Then
\begin{eqnarray*}
\lambda \phi =
(-F_3l_2+F_2l_3,-F_3l_1+F_1l_3,-F_2l_1+F_1l_2,-F_3l_3+F_2l_4,\\
-F_3l_2+F_1l_4,-F_2l_2+F_1l_3,F_1l_5-F_2l_6+F_3l_7)=(0,\dots,0)
\end{eqnarray*}

Since the matrix $B$ (the upper left block of $\phi$) has full rank
($=4$), we see that $l_1=l_2=l_3=l_4=0$.  Thus $l_5,l_6$ and $l_7$
must give the surjection in question.  But these satisfy the relation
$F_1l_5-F_2l_6+F_3l_7=0$, so the map $\O_X(-s-2)^{\oplus 3}\rightarrow
\L$ factors through the cokernel of
$(F_1,-F_2,F_3):\O_X(-2s-2)\rightarrow \O_X(-s-2)^{\oplus 3}$, and
this rank two cokernel (being non-split) has no invertible quotient.  Thus
there is no such surjection.  This concludes the proof of the proposition.
\end{proof}

\begin{bem}
The proof actually shows that we cannot get any structure of the
proposed type for higher dimensions as well.
\end{bem}

\begin{bem}  There is, however, a surjection onto a locally free sheaf
  of rank two (the cokernel of the map $(F_1,-F_2,F_3)$ occuring in
  the proof), showing that there is a thickening to a scheme of
  multiplicity {\em six}.
\end{bem}

\section{\label{kvotient-->multippelstruktur}Applications in
projective geometry: splittings of bundles and complete intersections}

Let $X \subset P$ be a regular embedding of smooth varieties, with normal
bundle $\N$ and Ideal $\I$.  Then $\N^{\vee}\iso \I/\I^2$.  Thus if we
have a quotient $\E$ of the conormal bundle we can define a new scheme $Y$
with Ideal $\J$ as follows:

\[\J:=\ker(\I\rightarrow \N^{\vee} \rightarrow \E).\]

\noindent Obviously, $\I^2 \subset \J \subset \I$, so $Y$ is contained in
$X^{(1)}$, the first infinitesimal neighbourhood of $X$ in $P$.  If
$\E$ is a bundle, then $Y$ is a Cohen-Macaulay scheme.  We will
primarily focus on the case $P=\PP^N$.\\

\subsection{Hartshorne's Conjecture}
A natural question is whether interesting properties of the bundle
$\E$ is reflected in the geometry of the scheme $Y$.  Recall the well
known conjectures

\begin{conj}[Hartshorne's Conjecture in codimension two]
Any smooth variety $X$ of codimension two in a projective space $\PP^N$
of dimension at least six is a complete intersection.
\end{conj}

\begin{conj}[The rank two bundle conjecture]
Any vector bundle of rank two on a projective space $\PP^N$ of
dimension at least six splits.
\end{conj}

Now consider a conjecture expressed entirely in terms of multiple
structures of degrees two and three:

\begin{conj}[Conjecture on triple linear schemes]
\label{tripconj}
Consider a linear subspace $\PP^N\subset \PP^{N+M}$ of dimension
$N\geq 6$ and (arbitrary) codimension $M$, and a degree three
Cohen-Macaulay scheme $Y$ satisfying 
\[\PP^N\subset Y\subset (\PP^N)^{(1)}.\]
Then there exists a degree two Cohen-Macaulay scheme $Z$ such that
\[\PP^N\subset Z \subset Y.\]
\end{conj}

\begin{thm}
\label{Hartshorneconj}
These three conjectures are equivalent for each $N$.
\end{thm}

\begin{proof}
The equivalence of Hartshorne's Conjecture and the Rank two Bundle
Conjecture is well known:  a smooth subvariety of codimension two in
$\PP^N, N\geq 6$ is subcanonical, which implies that its normal bundle
extends to the ambient space.  This bundle is split if and only if the
variety is a complete intersection.  Also, given any bundle of rank
two, sufficiently high twists of it will have smooth sections.\\

Now consider a setup as in the Conjecture on triple linear schemes.
There is a short exact sequence
\[0\rightarrow \I_Y/\I_X^2\rightarrow \I_X/\I_X^2\rightarrow
\E\rightarrow 0\]
where $X=\PP^N\subset \PP^{N+M}$.  According to Proposition
\ref{CMlokfri} $Y$ is Cohen-Macaulay if and only if $\E$ is locally
free.  By Theorem \ref{ekvlokfrilcm} there is a $1-1$ correspondence
between Cohen-Macaulay schemes $Z$ satisfying
\[X\subset Z\subset Y\]
and locally free sheaves $\L$ such that there is a surjection
\[\E\surj \L.\]
In this case, if both inclusions are strict, $Z$ has to be of degree
two and $\L$ has to be invertible.  Thus the kernel of $\E\rightarrow
\L$ is also invertible, so $\E$ splits as the sum of these two
bundles.  The $1-1$ correspondence says that we can also go the other
way; a splitting of $\E$ produces $Z$.  The only thing left to show is
that any bundle of rank two appears in such a sequence, or rather that
a twist of it does.\\

To see this, for a given $\E$, replace it by a
sufficiently high twist so that $\E(1)$ is generated by global sections.
Let $M=h^0(\E(1))$, and embed $X=\PP^N\subset\PP^{N+M}$ as a linear
subspace.  Then the conormal bundle of $X\subset\PP^{N+M}$
satisfies
\[\I_X/\I_X^2\iso \O_X(-1)^{\oplus M}\iso \O_X(-1)\otimes
H^0(\E(1))\surj \E\]
where the right hand map is the global section map twisted by $-1$.
Thus we get a sequence of the kind we need.
\end{proof}

\begin{bem}  We can also make a similar conjecture for some special
  multiple structures on a linear subspace of codimension two, which
  again will be equivalent to Hartshorne's Conjecture.  The degree,
  however, will be very much larger.  Explicitly, consider a vector
  bundle $\E$ on $\PP^N$, and assume that it is globally generated by
  its $M=h^0(\E)$ sections.  Restricting the conormal sheaf of the
  $(M-2)$nd infinitesimal neighbourhood of $\PP^N$ in $\PP^{N+2}$ to
  $\PP^N$, we get $\O_{\PP^N}(-(M-1))^{\oplus M}$, and $\E(-(M-1))$ is
  a quotient of this bundle.  We can now proceed as above to prove
  that Hartshorne's Conjecture is equivalent to the following
  conjecture:
\end{bem}

\begin{conj} Consider $\PP^N$, $N\geq 6$, as a linear subspace of
  codimension two in $\PP^{N+2}$.  Assume that $Y$ is a Cohen-Macaulay
  scheme satisfying
\[(\PP^N)^{(M)}\subset Y\subset (\PP^N)^{(M+1)}\]
such that the degree of $Y$ is equal to the degree of
$(\PP^N)^{(M)}+2$.  Then there is a Cohen-Macaulay scheme
$Z$ of degree equal to the degree of $(\PP^N)^{(M)}+1$ with
\[(\PP^N)^{(M)}\subset Z\subset Y.\]
\end{conj} 

\begin{bem}
These two conjectures on multiple structures consider ``extreme''
invariants:  low degree and low codimension.  The following more
general conjecture can be proven to be equivalent to Hartshorne's
Conjecture in exactly the same way:
\end{bem}

\begin{conj}  Consider a linear subspace $X=\PP^N\subset \PP^{N+M}$ of
  dimension $N\geq 6$.  Let $W$ be a Cohen-Macaulay multiple monomial
  scheme with $W_{\red}=X$, and let $Y$ be a Cohen-Macaulay multiple
  scheme satisfying
\[\I_W\I_X\subset \I_Y\subset \I_W\]
and $\deg Y=\deg W +2$.  Then there exists a Cohen-Macaulay scheme
$Z$, of degree equal to $\deg W +1$, such that
\[W\subset Z\subset Y.\]
\end{conj}

\begin{bem}
Any counter-example to any of these statements for $N\leq 5$ also
produces a non-split bundle of rank two.
\end{bem}

\begin{bem}
It is well known that any bundle of rank $r$ on $\PP^N$ is a quotient
of a sum of $N+r$ line bundles.  Thus in order to {\em prove}
Hartshorne's conjecture, it would be sufficient to classify triple
Cohen-Macaulay schemes on a linear subspace $X$ of dimension six in
$\PP^14$, contained in the first infinitesimal neighbourhood.  If
Hartshorne's conjecture holds, the only ones are of the following
form:  

If a given rank two bundle splits, the sequence
defining its triple structure would be 

\[ 0\rightarrow \I/\I_X^2\rightarrow \I_X/\I_X^2\iso \O(-1)^{\oplus 8}
\rightarrow \O(a)\oplus \O(b)\rightarrow 0\]
where the right hand map is given by a full-rank matrix 

\[\left( \begin{array}{cccccccc} f_0 & f_1 & f_2 & f_3 & f_4 &f_5 &f_6 &f_7\\
                       g_0 & g_1 & g_2 & g_3 & g_4 &g_5 &g_6 &g_7
                       \end{array} \right)\]

Let $I_X=(x_0,\cdots, x_7)$ and $M_{ij}=f_ig_j-f_jg_i$ for all $i>j$.
Then the ideal $I$ can be computed from the above as

\[I=(x_kM_{ij}-x_jM_{ik}+x_iM_{jk})+I_X^2, \mbox{ } 7\geq i>j>k \geq 0.\]

\end{bem}

Because of the link with bundles of rank two, we are interested in
schemes of degree three.  For reduced,
equidimensional schemes there is a short list of examples.  One might
wonder how these possible schemes of degree three compare with the
reduced schemes; this question is considered in Section \ref{reddegtre}.\\

The last part of the proof of Theorem \ref{Hartshorneconj} can be generalized to a much wider setting,
though the small degree part is removed (as well as the condition that
$\E(1)$ is generated by global sections):

\begin{prop}
  Let $\E$ be a bundle on a projective variety $X$.  Then there is a
  projective embedding $i:X\hookrightarrow \PP^N$ with normal bundle
  $\N_i$, and a surjection $\N_i^{\vee}\surj \E$.
\end{prop}
 \begin{proof}
  Let $\pi:\PP(\E^{\vee} \oplus \O_X)\rightarrow X$ be the projective
  completion of the total space of the bundle $\E^{\vee}$.  The normal bundle of the zero
  section $s$, $\N_s$, is isomorphic to $\E^{\vee}$.  The projective completion
  is again a 
  projective variety, and can thus be embedded $j:\PP(\E^{\vee} \oplus
  \O_X)\hookrightarrow \PP^N$.  Define $i=j \circ s$.  We have a short
  exact sequence (see \cite{SGA6}, Proposition VII.1.7.  The $\N$ there is our
  $\N^{\vee}$.)  
\[0 \rightarrow s^{\ast}\N_j^{\vee} \rightarrow \N_i^{\vee}
  \rightarrow \N_s^{\vee} \rightarrow 0.\]

But $\E\iso \N_s^{\vee}$.
\end{proof}

There is a way of doing this more directly for an already (regularly)
embedded 
variety $i:X\hookrightarrow \PP^N$, but this again requires $\E(1)$ to
be globally generated.  Let $\N$ be the normal bundle of $i$.

\begin{prop}
\label{buntkvotient}
If $\E$ is a bundle on $X$ such that $\E(1)$ is generated by global
sections, then there is a further embedding $X
\stackrel{i}{\subseteq}\PP^N \stackrel{j}{\subset}\PP^{N+M}$, where $j$ is a
linear embedding, such that $\E$ is a quotient of the conormal bundle of
$j\circ i$.
\end{prop}
\begin{proof}
Note first that the normal bundle of $\PP^N$ in $\PP^{N+M}$ is isomorphic
to $\O(1)^{\oplus M}$.  We have an exact sequence
\[0\rightarrow \N \rightarrow \N_{X/\PP^{N+M}}\rightarrow
i^{\ast}\N_{\PP^N/\PP^{N+M}} \rightarrow 0.\]
We want to show that this sequence splits.  Using induction, we can
assume $M=1$.  Since the normal bundle depends only on a neighbourhood of the
embedded variety, and because $\N_{\PP^N/\PP^{N+1}}\iso\O(1)$, this
exact sequence is equivalent to
\[0\rightarrow \N \rightarrow \N_{X/\PP^{N+1}\setminus \{P\}}\rightarrow
\O_X(1) \rightarrow 0\]
where $P\in \PP^{N+1}\setminus \PP^N$ is a point.  We thus have a
retraction $p:\PP^{N+1}\setminus \{P\} \rightarrow 
\PP^N$.  Since $p\circ j \circ i=i$, a regular embedding, we have the
sequence (see \cite{Ful98} B 7.5.)
\[0\rightarrow (j\circ i)^{\ast}\T_{\PP^{N+1}\setminus \{P\}/\PP^N}
\rightarrow \N_{X/\PP^{N+1}\setminus \{P\}}\rightarrow \N\rightarrow 0.
\]
This is easily seen to give a splitting of the sequence above.  Thus 
\[\N_{X/\PP^{N+M}}\iso \N \oplus \O_X(1)^{\oplus M}.\]  
Choose $M=h^0(\E(1))$, so that we have a morphism 
\[\O_X^{\oplus M} \iso\H^0(\E(1))\otimes \O_X
\stackrel{\alpha(1)}{\surj} \E(1)\] 
and thus
\[\N_{X/\PP^{N+M}}^{\vee}\iso \N^{\vee} \oplus
  \O_X(-1)^{\oplus M}\stackrel{(0,\alpha)}{\surj} \E\]
as claimed.
\end{proof}

\begin{defn}[The multiple structure associated to a bundle]
Let us first consider the case $X=\PP^N$, with a bundle $\E$ such that
$\E(1)$ is globally generated.  
The following multiple structure $Y$ will be called {\em the multiple
  structure associated to $\E$}:  Embed $\PP^N\subset \PP^{N+M}$ as a
  linear subspace, $M=h^0(\E(1))$.  Then $\E$ is a quotient of the
  conormal bundle $\N^{\vee}\iso \I/\I^2$.  The Ideal $\J$ defined by 
\[\J:=\ker(\I\surj\I/\I^2 \surj \E)\]
defines a Cohen-Macaulay scheme $Y$, with $Y_{\red}=\PP^N$.  This is
the scheme constructed, for $\E$ of rank two, in the proof of Theorem
\ref{Hartshorneconj}.\\

In general, for an embedded scheme $X\subset\PP^N$, and $\E$ as above, we
define $\J$ and $Y$ in the same manner, but using the construction
from Proposition \ref{buntkvotient}.  Note that this reduces to the
above in case $X=\PP^N$.
\end{defn}

This notion behaves well under the operation direct sum, and also
under restriction if some extra conditions are met.  We assume
throughout that $\E(1)$ is a globally generated bundle on $X\subset
\PP^N$.  Thus the multiple scheme associated to $\E$ is defined.

\begin{thm}
\label{koproduktteoremet}
Suppose $\E$ splits as a sum of two bundles $\E\iso \E_1\oplus \E_2$.
Then the embedding and the multiple scheme associated to $\E$ are the
cofibered coproduct of the embeddings and multiple schemes associated
to $\E_1$ and $\E_2$.
\end{thm}
\begin{proof}
Define $M_i=h^0(\E_i(1)),M=M_1+M_2=h^0(\E(1))$.  There is a diagram of
inclusions

\[\xymatrix{
X\ar[rd]\ar[ddd]\ar[rrr]&                        &  & Y_1\ar[ddd]\ar[dl]\\ 
        &\PP^N\ar[r]\ar[d]\ar[rd]&\PP^{N+M_1}\ar[d] & \\
        &\PP^{N+M_2}\ar[r]       &\PP^{N+M}&\\
Y_2 \ar[rrr]\ar[ur]&&&Y\ar[ul]}\]

\noindent Obviously, the inner square is Cartesian, as well as the upper and
left squares.  To see that the right square is Cartesian, we need to
show that $Y\bigcap \PP^{N+M_1}=Y_1$ as subschemes of $\PP^{N+M}$.
Let $\I,\J,\J_1$ be the Ideals of $X,Y,Y_1$ in $\O_{\PP^{N+M}}$.  We
need to show that $J_1=J+I_{\PP^{N+M_1}/\PP^{N+M_1+M_2}}$, and since
both sides contain $I^2$, it will suffice to show equality modulo
this ideal.  By
the proof of Proposition \ref{buntkvotient} $\I/\I^2\iso
\N^{\vee}\oplus \O_X(-1)^{\oplus M_1}\oplus \O_X(-1)^{\oplus M_2}$.  Then
we have a short exact sequence

\[0\rightarrow \J_1/\I^2\rightarrow \I/\I^2\rightarrow \E_1\rightarrow 0,\]

or, using the splitting of $\I/\I^2$,
\[0\rightarrow \J_1/\I^2\rightarrow \N^{\vee}\oplus \O_X(-1)^{\oplus
  M_1}\oplus \O_X(-1)^{\oplus M_2} \rightarrow \E_1\rightarrow 0.\]

Since the right hand map splits as the sum of two zero maps and a
surjection $\alpha_1:\O_X(-1)^{\oplus M_1}\rightarrow \E_1$, we get an
isomorphism 
\[\J_1/\I^2\iso \N^{\vee}\oplus \ker \alpha_1\oplus \O_X(-1)^{\oplus
  M_2}.\]

Likewise 

\[\J/\I^2\iso \N^{\vee}\oplus \ker \alpha_1\oplus \ker \alpha_2.\]

Thus we see that the difference in the generating sets of the ideals
$J$ and $J_1$ is just the linear polynomials in the $\O_X(-1)^{\oplus
  M_2}$ part, which are the generators for the ideal
$I_{\PP^{N+M_1}/\PP^{N+M_1+M_2}}$.  The proof that the right hand
square is Cartesian is complete.\\

The lower square is thus Cartesian by symmetry.  Since all the small
  squares are Cartesian, the outer square must be Cartesian as well.
  The proof is complete. 
\end{proof}

\begin{thm}
Let $\E$ be a bundle on $Z$, $X\subset Z\subset \PP^N$ smooth schemes
and regular embeddings.  Assume $\E(1)$ is globally generated on $Z$,
and that the canonical map $H^0(\E(1))\rightarrow H^0(\E(1)|_X)$ is
surjective (e.g. if $\H^1(\E(1)\otimes \I_{X/Z})=0$).  Then the multiple scheme
associated to $\E|_X$ is the intersection of the multiple scheme
associated with $\E$ and the cone over $X$ (in the embedding of
$\PP^N\subset \PP^{N+M}$ associated to $\E$).  This intersection is
contained in a linear subspace of
codimension $h^0(\E(1)\otimes \I_{X/Z})$.
\end{thm}
\begin{proof}
Note first that the conditions imply that
\[0\rightarrow \H^0(\E(1)\otimes \I_{X/Z})\rightarrow \H^0(\E(1))
\rightarrow \H^0(\E(1)|_X)\rightarrow 0\]
is exact.\\

When taking the multiple scheme associated with $\E(1)$ we
choose a basis for $\H^0(\E(1))$ and use this basis as a new set of
variables, generating the ideal of $\PP^N$ in $\PP^{N+M}$, where
$M=h^0(\E(1))$.  Let $\I,\J$ be the Ideals of $X,Z$ in $\PP^{N+M}$.
Choose a vector space splitting 
$\H^0(\E(1))\iso\H^0(\E(1)\otimes \I_{X/Z}) \oplus \H^0(\E(1)|_X)$.
Consider the two sequences 
\[0\rightarrow \I'/\I^2\rightarrow \N_{Z/\PP^N}^{\vee}\oplus \H^0(\E(1))\otimes
\O_Z(-1)\rightarrow \E\rightarrow 0\]
and
\[0\rightarrow \J'/\J^2\rightarrow \N_{X/\PP^N}^{\vee}\oplus
(\H^0(\E(1)_{|X})\oplus \H^0(\E(1)\otimes \I_{X/Z}))\otimes
\O_X(-1)\rightarrow \E|_X\rightarrow 0.\]

Clearly the kernel of the lower sequence is generated by the conormal
bundle of $X$ in $\PP^N$ (corresponding to the ideal of the
cone over $X$) and the kernel of the upper sequence (intersecting with
the scheme associated to $\E$).  The term $\H^0(\E(1)\otimes \I_{X/Z})\otimes
\O_X(-1)$ in the lower sequence tells us that there are
$h^0(\E(1)\otimes \I_{X/Z})$ linear polynomials in the ideal.  The
remaining equations are the equations of the scheme associated to
$\E|_X$ in the {\em embedding} associated to $\E|_X$.

\end{proof}

\subsection{Extensions and splittings of bundles; the first
  infinitesimal neighbourhood}

In this section we will generalize the discussion about
Hartshorne's Conjecture to include a much wider class of questions
concerning the splittings and extensions of bundles.  We will also
give some examples, including equations for the scheme
associated to the Horrocks-Mumford bundle.\\

We let $\E$ be a bundle on $\PP^N$, of rank $r$, such that $\E(1)$ is
globally generated.  Define $M=h^0(\E(1))$ and let $Y$ be the scheme
associated to $\E$; there is a sequence
\[0\rightarrow \I_Y/\I_X^2\rightarrow \I_X/\I_X^2\rightarrow
\E\rightarrow 0\]
where all Ideals are Ideals in $\O_{\PP^{N+M}}$.
\begin{thm}
Suppose $N\geq 2$.  The bundle $\E$ splits as a sum of line bundles if
and only if there is a sequence of subschemes
\[X=\PP^N\subset Z_1\subset Z_2\subset \dots\subset Z_{r-1}\subset Y\]
where each $Z_i$ is a Cohen-Macaulay scheme, and the degree of $Z_i$ is
$i+1$. 
\end{thm}
\begin{proof}
For each $Z_i$ define the locally free $\O_X$-Module $\F_i$, of rank
$i$, by
\[0\rightarrow \I_{Z_i}/\I_X^2\rightarrow \I_X/\I_X^2\rightarrow
\F_i\rightarrow 0.\]
By Theorem \ref{ekvlokfrilcm}, there is a sequence of surjective maps
(or a cofiltration of $\E$)
\[\E\surj\F_r\surj\F_{r-1}\surj\dots \surj \F_1\]
where the kernel of each map is an invertible sheaf, 
\[0\rightarrow \M_i \rightarrow \F_i\rightarrow \F_{i-1}\rightarrow 0.\]

Since $N\geq 2$, all $\Ext^1_{\O_{\PP^N}}(\G_1,\G_2)=0$ for invertible
sheaves $\G_1,\G_2$.  Especially, the sequence
\[0\rightarrow \M_2\rightarrow \F_2\rightarrow \F_1\rightarrow 0\]
splits.  Inductively we get
\[\F_i\iso \oplus_{j \leq i}\M_j\]
(where $M_1:=\F_1$) from the exact sequence
\[0\rightarrow \M_i\rightarrow \F_i\rightarrow \F_{i-1}\rightarrow 0\]
and
\[\Ext^1(\F_{i-1},\M_i)\iso\Ext^1(\oplus_{j \leq i-1}\M_j,\M_i)\iso
\oplus_{j \leq i-1}\Ext^1(\M_j,\M_i)=0.\]

Thus
\[\E\iso \oplus \M_i.\]

The opposite implication can be shown by reversing this argument,
again using Theorem \ref{ekvlokfrilcm} to conclude.
\end{proof}

One of the difficulties inherent in this theory is that the
codimension must be rather large if we want interesting examples.  For
instance 

\begin{prop}\label{lavcodim}
Let $\PP^N$ be linearly embedded in $\PP^{N+M}$, $M\leq N$, and let
there be given a Cohen-Macaulay nilpotent scheme $Y$ with support equal to this
$\PP^N$, contained in the first infinitesimal neighbourhood.  Then $Y$
has a Cohen-Macaulay subscheme with the same support, of degree one less.
\end{prop}
\begin{proof}
This is a consequence of a result of Faltings' \cite{Fal81}, which
states (among other things) that a surjection from a free
$\O_{\PP^N}$-Module of rank $\leq N$ onto a locally free Module is split. 
\end{proof}

Thus the first interesting case is $\PP^2\subset \PP^5$, since all
bundles on $\PP^1$ split.

\begin{example}
Consider the embedding of $X=\PP^2 \subset \PP^5$.  To find interesting
examples of Cohen-Macaulay multiple structures of degree three, we
need bundles of rank
two that are quotients of $\O_X(-1)^{\oplus 3}$.  We start with the dual
situation; more precisely we have a short exact sequence

\[ 0 \rightarrow \Q \rightarrow \O_X(1)^{\oplus 3} \rightarrow \O_X(a+2)
\rightarrow 0,\]
\noindent
 where $a \geq 0$ and the last map is given by three homogeneous
 polynomials $f_1,f_2,f_3$ of degree $a+1$ in the variables
 $z_0,z_1,z_2$.  The surjectivity of the map is equivalent to the
 condition $\bigcap V_+(f_i)=\emptyset$.  $\Q$ is the kernel of this map,
 a bundle of rank two.  We dualize the above sequence to get 

\[0 \rightarrow \O_X(-a-2) \rightarrow \O_X(-1)^{\oplus 3} \rightarrow
 \F \rightarrow 0, \]
\noindent
$\F=\Q^{\vee}$.  Now the first map is given by $1 \mapsto
(f_1e_1,f_2e_2,f_3e_3)^t$, where $1$ is a basis for $\O_X(-a-2)$ and the
$e_i$ are a basis for $\O_X(-1)^{\oplus 3}$.  \\

  In our situation, we fit this into our standard sequence 

\[ 0 \rightarrow \I_Y/\I_X^2 \rightarrow \I_X/\I_X^2 \rightarrow \F
\rightarrow 0. \]
\noindent
The basis element $e_i$ corresponds to the monomial $x_i$ in $I_X$.
Obviously, $\F(1)$ is globally generated by three sections, whereas $\F$
is not globally generated (in fact, $\F$ doesn't have any global
sections).  The ideal of $Y$ is given by the equations

\[I_Y=(f_1x_1+f_2x_2+f_3x_3,(x_1,x_2,x_3)^2). \]

Compare this example with the multiplicity five structure from
Proposition \ref{nystruktur}.
\end{example}

\begin{bem}  Note that there is no other Cohen-Macaulay triple structure on
  $\PP^2\subset \PP^5$ that is contained in $(\PP^2)^{(1)}$, without
  a double Cohen-Macaulay substructure.
\end{bem}

\begin{example}[The split case]
This example will shed light on Proposition \ref{lavcodim}, Theorem
\ref{koproduktteoremet} and Conjecture \ref{tripconj}.\\

Consider the following setup:\\

Let $X=\PP^n \subset \PP^{n+m}$ be the linear embedding associated to
$\E \iso \O_X(a) \oplus \O_X(b), a,b\geq 0$ by the above method (the case
$a=-1$ being simpler).  In particular, $m=h^0(\E(1))$.  Then the
corresponding multiplicity 3 structure has equations 
\[(x_{I_1}z^{I_2} - x_{I_3}z^{I_4},x_{J_1}z^{J_2} -
x_{J_3}z^{J_4},(x_I,x_J)^2)| I_1+I_2=I_3+I_4,J_1+J_2=J_3+J_4\]
where $X$ has coordinates $z_0,...,z_n$, $z^I:=z_0^{i_0}...z_n^{i_n}$
for the $n+1$-tuple $I=(i_0,...,i_n)$ of weight $\|I\|=\Sigma i_l
=a+1$, $z^J:=z_0^{j_0}...z_n^{j_n}$ for the $n+1$-tuple $I=(j_0,...,j_n)$ of weight $\|J\|=\Sigma j_l=b+1$
and the ideal $I$ of $X$ is generated by the $x_I$ and the $x_J$.
Each tuple with letter $I$ has weight $a+1$, and similarly for $J$.
In the equation $I_1+I_2=I_3+I_4$, all operations are operations on
{\em tuples}, so the two expressions are equal term by term.
\begin{proof}
Recall the morphism 
\[H^0(\E(1))\otimes \O_X \rightarrow \E(1)\]
which, twisted with $-1$, gives the exact sequence
\[ \xymatrix{
0 \ar[r]& \J/\I^2 \ar[r]& \I/\I^2 \iso \O_X(-1)^{\oplus m} \ar[r]^<<<<<{\alpha} &
\E \ar[r] & 0 }. \]
$\J$ is the Ideal we are interested in.  $\alpha$ consists of two direct
summands:  the morphism defined
by a generating set of the vector space $H^0(\O_X(a+1))$, which we
take to be the monomials of degree $a+1$ in the $z_i$'s, and the other
morphism similarly defined by the monomials of degree $b+1$. 
These parts contribute with all
relations on the monomials of degree $a+1$ and of degree $b+1$.  For
each monomial $z^I:=z_0^{i_0}...z_n^{i_n}$, let $x_I$ be the
corresponding variable in the homogeneous coordinate ring of
$\PP^{n+m}$.  We want to compute the relations
among the $x_I$ with  coefficients in the $z_i$.  Compare with the
$(a+1)$-uple embedding associated to the complete  linear system $\|
\O_X(a+1) \| $, given on coordinate rings by
\[\begin{array}{ccc}
k[T_I \| \mbox{I is an n-tuple of total weight a+1}]& \rightarrow & k[z_0,...,z_n].
\\
&T_I \mapsto z^I &

\end{array} \]
 The kernel of this map, i.e. the ideal of the $(a+1)$-uple embedding, is
 generated by the elements $T_{I_1}T_{I_2}-T_{I_3}T_{I_4},
 I_1+I_2=I_3+I_4$.  Translating this back to the $x_I$, the relations
 among them with coefficients in the $z_i$ are exactly $x_{I_1}z^{I_2}
 - x_{I_3}z^{I_4},I_1+I_2=I_3+I_4$.  Note that also the relations with
 $I_1=I_3, I_2=I_4$ must be included, and that the two relations
 $x_{I_1}z^{I_2} - x_{I_3}z^{I_4}$ and $x_{I_2}z^{I_1}-x_{I_4}z^{I_3}$
 are different.  Thus the number of equations is larger than the
 number of equations for the $(a+1)$-uple embedding.  Giving the other
 direct summand, corresponding to monomials of degree $b+1$, the same
 treatment proves the statement. 
\end{proof}

Since the bundle in questions splits, we expect to find a
Cohen-Macaulay subscheme of degree two.  There are two choices for
this subschemes, one for each summand.  One possibility has ideal

\[(x_{I_1}z^{I_2} - x_{I_3}z^{I_4},x_J,(x_I)^2)| I_1+I_2=I_3+I_4.\]

Switching the role of $I$ and $J$ produces the other.
\end{example}

\begin{bem}  By Hartshorne's Conjecture we expect that the above is
  essentially the only example of triple structures in the first
  infinitesimal neighbourhood of a linear subspace of dimension at
  least six.
\end{bem}

\begin{example}
\label{HorrMum}
The famous Horrocks-Mumford bundle on $\PP^4=X$ is particularly
interesting because it is essentially the only known example of a
non-split rank 2 bundle on $\PP^4$ in characteristic zero.  So let $\E$ be this bundle,
twisted so that $\E(1)$ is generated by global sections, whereas $\E$
is not.  Then the minimal free
  resolution of $\E$ is given by Manolache \cite{Man88},\cite{Man89}
  (including explicit matrices) as

 \[ 0 \rightarrow \O_X(-5)^{\oplus 2} \rightarrow \O_X(-3)^{\oplus 20}
  \rightarrow \O_X(-2)^{\oplus 35} \rightarrow \O_X(-1)^{\oplus 15} \oplus
  \O_X^{\oplus 4} \rightarrow \E \rightarrow 0.\]
\noindent
 Note that $\dim \Gamma (\E(1))
  =35$, so that we will consider $X=\PP^4 \subset \PP^{39}$. Consider
  the composite
\[ \xymatrix{
\O_X(-1)^{\oplus 35} \ar[r]^>>>>{\phi} & \O_X(-1)^{\oplus 15} \oplus \O_X^{\oplus
  4} \ar[r] & \E }.\]
$\phi$ is the direct sum of two maps:  One is an isomorphism
  $\O_X(-1)^{\oplus 15} \iso \O_X(-1)^{\oplus 15}$.  The other is four
  copies of the right hand map in the Euler sequence

 \[ 0 \rightarrow \Omega_X \rightarrow \O_X(-1)^{\oplus 5} \rightarrow \O_X
  \rightarrow 0. \]
  
\noindent
Using the computer system Macaulay2, \cite{Mac2}, we find the equations defining the
kernel $\J/\I_X^2$ from the sequence
\[0\rightarrow \J/\I_X^2\rightarrow \O_X(-1)^{\oplus 35} \rightarrow \E
\rightarrow 0.\]
Together with the equations $(x_1,...,x_{35})^2$
they generate the ideal corresponding to $\J$.  The remaining equations are

\[ \begin{array}{cc} 
-{z_4} {{x}}_{{19}}+{z_3} {{x}}_{{20}} &
-{z_4} {{x}}_{{24}}+{z_3} {{x}}_{{25}}\\
-{z_4} {{x}}_{{29}}+{z_3} {{x}}_{{30}} &
-{z_4} {{x}}_{{34}}+{z_3} {{x}}_{{35}}\\
-{z_4} {{x}}_{{5}}+{z_2} {{x}}_{{17}}+{z_3} {{x}}_{{21}}&
{z_3} {{x}}_{{2}}-{z_4} {{x}}_{{11}}+{z_2} {{x}}_{{18}}\\
{z_2} {{x}}_{{3}}-{z_4} {{x}}_{{16}}-{z_3} {{x}}_{{22}}&
-{z_3} {{x}}_{1}+{z_2} {{x}}_{{12}}-{z_4} {{x}}_{{20}}\\
-{z_3} {{x}}_{{18}}+{z_2} {{x}}_{{19}}&
-{z_3} {{x}}_{{23}}+{z_2} {{x}}_{{24}}\\
-{z_3} {{x}}_{{28}}+{z_2} {{x}}_{{29}}&
-{z_3} {{x}}_{{33}}+{z_2} {{x}}_{{34}}\\
-{z_4} {{x}}_{{18}}+{z_2} {{x}}_{{20}}&
-{z_4} {{x}}_{{23}}+{z_2} {{x}}_{{25}}\\
-{z_4} {{x}}_{{28}}+{z_2} {{x}}_{{30}}&
-{z_4} {{x}}_{{33}}+{z_2} {{x}}_{{35}}\\
-{z_3} {{x}}_{{4}}+{z_1} {{x}}_{{16}}+{z_4} {{x}}_{{23}}&
{z_3} {{x}}_{{9}}+{z_4} {{x}}_{{13}}+{z_2} {{x}}_{{15}}+{z_1} {{x}}_{{21}}-{z_4}
{{x}}_{{28}}\\
{z_4} {{x}}_{{4}}-{z_2} {{x}}_{{6}}+{z_1} {{x}}_{{22}}&
{z_2} {{x}}_{1}-{z_3} {{x}}_{{15}}+{z_1} {{x}}_{{17}}\\
{z_1} {{x}}_{{2}}-{z_4} {{x}}_{{19}}-{z_2} {{x}}_{{21}}&
{z_4} {{x}}_{{9}}+{z_3} {{x}}_{{10}}+{z_1} {{x}}_{{12}}+{z_2} {{x}}_{{16}}-{z_4} {{x}}_{{34}}\\
{z_2} {{x}}_{{7}}+{z_1} {{x}}_{{8}}+{z_4} {{x}}_{{15}}+{z_3} {{x}}_{{16}}-{z_2} {{x}}_{{32}}&
-{z_4} {{x}}_{1}+{z_1} {{x}}_{{9}}-{z_2} {{x}}_{{23}}\\
{z_2} {{x}}_{{8}}+{z_3} {{x}}_{{12}}+{z_1} {{x}}_{{14}}+{z_4} {{x}}_{{21}}-{z_3} {{x}}_{{27}}&
{z_1} {{x}}_{1}-{z_4} {{x}}_{{8}}+{z_3} {{x}}_{{24}}
\end{array} \]

\[\begin{array}{cc}
-{z_2} {{x}}_{{5}}+{z_1} {{x}}_{{11}}-{z_3} {{x}}_{{19}}&
{z_1} {{x}}_{{3}}-{z_3} {{x}}_{{6}}+{z_4} {{x}}_{{25}}\\
-{z_4} {{x}}_{{5}}+{z_1} {{x}}_{{18}}+{z_3} {{x}}_{{21}}&
-{z_2} {{x}}_{{22}}+{z_1} {{x}}_{{23}}\\
-{z_2} {{x}}_{{27}}+{z_1} {{x}}_{{28}}&
-{z_2} {{x}}_{{32}}+{z_1} {{x}}_{{33}}\\
-{z_3} {{x}}_{{17}}+{z_1} {{x}}_{{19}}&
-{z_3} {{x}}_{{22}}+{z_1} {{x}}_{{24}}\\
-{z_3} {{x}}_{{27}}+{z_1} {{x}}_{{29}}&
-{z_3} {{x}}_{{32}}+{z_1} {{x}}_{{34}}\\
-{z_4} {{x}}_{{17}}+{z_1} {{x}}_{{20}}&
-{z_4} {{x}}_{{22}}+{z_1} {{x}}_{{25}}\\
-{z_4} {{x}}_{{27}}+{z_1} {{x}}_{{30}}&
-{z_4} {{x}}_{{32}}+{z_1} {{x}}_{{35}}\\
{z_3} {{x}}_{{3}}-{z_1} {{x}}_{{10}}+{z_0} {{x}}_{{21}}&
-{z_3} {{x}}_{{3}}+{z_2} {{x}}_{{4}}-{z_4} {{x}}_{{7}}+{z_1} {{x}}_{{10}}\\
{z_1} {{x}}_{{5}}-{z_2} {{x}}_{{14}}+{z_0} {{x}}_{{16}}&
-{z_4} {{x}}_{{2}}+{z_1} {{x}}_{{5}}+{z_3} {{x}}_{{13}}-{z_2} {{x}}_{{14}}\\
{z_1} {{x}}_{{6}}+{z_0} {{x}}_{{7}}+{z_3} {{x}}_{{14}}+{z_4} {{x}}_{{18}}-{z_1} {{x}}_{{31}}&
{z_4} {{x}}_{{4}}-{z_3} {{x}}_{{5}}-{z_2} {{x}}_{{6}}+{z_0} {{x}}_{{8}}\\
-{z_2} {{x}}_{1}+{z_0} {{x}}_{{3}}-{z_4} {{x}}_{{14}}+{z_3} {{x}}_{{15}}&
{z_1} {{x}}_{{7}}+{z_2} {{x}}_{{11}}+{z_0} {{x}}_{{13}}+{z_4} {{x}}_{{24}}-{z_2} {{x}}_{{26}}\\
{z_4} {{x}}_{{10}}+{z_3} {{x}}_{{11}}+{z_0} {{x}}_{{14}}+{z_2} {{x}}_{{22}}-{z_3} {{x}}_{{26}}&
{z_0} {{x}}_{{5}}-{z_3} {{x}}_{{7}}+{z_2} {{x}}_{{23}}\\
{z_3} {{x}}_{{2}}-{z_1} {{x}}_{{4}}-{z_4} {{x}}_{{11}}+{z_0} {{x}}_{{15}}&
{z_4} {{x}}_{{6}}+{z_0} {{x}}_{{10}}+{z_2} {{x}}_{{13}}+{z_3} {{x}}_{{17}}-{z_4} {{x}}_{{31}}\\
{z_0} {{x}}_{1}-{z_3} {{x}}_{{18}}-{z_4} {{x}}_{{22}}&
{z_0} {{x}}_{{6}}+{z_4} {{x}}_{{12}}+{z_1} {{x}}_{{15}}+{z_3} {{x}}_{{23}}-{z_4} {{x}}_{{27}}\\
{z_3} {{x}}_{{8}}+{z_2} {{x}}_{{9}}+{z_0} {{x}}_{{11}}+{z_4} {{x}}_{{17}}-{z_3} {{x}}_{{33}}&
{z_0} {{x}}_{{2}}-{z_2} {{x}}_{{10}}+{z_3} {{x}}_{{24}}\\
-{z_4} {{x}}_{{3}}+{z_0} {{x}}_{{12}}-{z_3} {{x}}_{{19}}&
-{z_2} {{x}}_{{2}}+{z_0} {{x}}_{{9}}-{z_4} {{x}}_{{25}}\\
{z_0} {{x}}_{{4}}-{z_1} {{x}}_{{13}}+{z_4} {{x}}_{{20}}&
-{z_3} {{x}}_{{4}}+{z_0} {{x}}_{{17}}+{z_4} {{x}}_{{23}}\\
{z_3} {{x}}_{{9}}+{z_4} {{x}}_{{13}}+{z_2} {{x}}_{{15}}+{z_0} {{x}}_{{22}}-{z_4} {{x}}_{{28}}&
-{z_1} {{x}}_{{26}}+{z_0} {{x}}_{{27}}\\
-{z_1} {{x}}_{{31}}+{z_0} {{x}}_{{32}}&
-{z_2} {{x}}_{{16}}+{z_0} {{x}}_{{18}}\\
-{z_2} {{x}}_{{21}}+{z_0} {{x}}_{{23}}&
-{z_2} {{x}}_{{26}}+{z_0} {{x}}_{{28}}\\
-{z_2} {{x}}_{{31}}+{z_0} {{x}}_{{33}}&
-{z_3} {{x}}_{{16}}+{z_0} {{x}}_{{19}}\\
-{z_3} {{x}}_{{21}}+{z_0} {{x}}_{{24}}&
-{z_3} {{x}}_{{26}}+{z_0} {{x}}_{{29}}\\
-{z_3} {{x}}_{{31}}+{z_0} {{x}}_{{34}}&
-{z_4} {{x}}_{{16}}+{z_0} {{x}}_{{20}}\\
-{z_4} {{x}}_{{21}}+{z_0} {{x}}_{{25}}&
-{z_4} {{x}}_{{26}}+{z_0} {{x}}_{{30}}\\
-{z_4} {{x}}_{{31}}+{z_0} {{x}}_{{35}} &
\end{array}
\]

\noindent
Thus we have a large number of equations (705), and the codimension is 35.  

\end{example}

\subsection{Reduced schemes of degree three}
\label{reddegtre}
In this section we will consider reduced, equidimensional projective
schemes of low degree.  For degree one the only such variety is linear
space itself.  For degree two there are two possibilities: a quadric
hypersurface (which may or may not be smooth) and a union of two
linear subspaces, which may intersect in any given dimension.  For
degree three the classification was done in the anonymous note
\cite{XXX}.  Roughly speaking, there is one irreducible variety of
degree three that is not a hypersurface, namely the Segre embedding of
$\PP^1\times \PP^2 \subset \PP^5$, together with (repeated) hyperplane
sections and cones over these.  For degree four Swinnerton-Dyer
classified the irreducible varieties in \cite{Swi73}.  For
degrees five to eight all varieties have been classified in a series of
papers by Ionescu, see \cite{Ion90} and the references therein.  Because of
the link with Hartshorne's Conjecture (from Conjecture \ref{tripconj} and
Theorem \ref{Hartshorneconj}) we will focus our attention on
the degree three case.  A complete list will be given explicitly only
in dimension four, and the Hilbert polynomials will be compared to the
Hilbert polynomial of the scheme associated to the Horrocks-Mumford
bundle.\\

\begin{prop}
\label{degreethree}
 The irreducible, nondegenerate varieties of degree three
  are the following:
\begin{itemize}
\item[(i)] A cubic hypersurface in $\PP^n$.
\item[(ii)] The Segre variety $\PP^1\times \PP^2 \subset \PP^5$.
\item[(iii)] An irreducible hyperplane section of the Segre variety in
  (ii).
\item[(iv)] A twisted cubic in $\PP^3$; this is an irreducible
  hyperplane section of the variety in (iii).
\item[(v)] A cone over one of the varieties in (ii)-(iv).
\end{itemize}
\end{prop}

Of course, a cone over a (cubic) hypersurface is again a (cubic)
hypersurface.  \\

The Hilbert polynomials of these varieties are most conveniently expressed
in terms of the Hilbert polynomials of projective spaces.

\begin{prop}Let $P_n(d):=\left(\begin{array}{c}n+d \\
      d\end{array}\right)$ be the Hilbert polynomial of projective
      $n$-space.  Then the Hilbert polynomials of (i)-(v) in Proposition
      \ref{degreethree} are
\begin{itemize}
\item[(i)] $3P_{n-1}-3P_{n-2}+P_{n-3}$
\item[(ii)] $3P_3-2P_2$
\item[(iii)] $3P_2-2P_1$
\item[(iv)] $3P_1-2P_0$
\item[(v)] If the dimension of the cone is $n$, the Hilbert polynomial
  is $3P_n-2P_{n-1}$
\end{itemize}
\end{prop}

\begin{proof}
For the first two, use the resolution of the structure sheaf.  Then (iii)
and (iv) follow from (ii) since the resolutions come from the resolution of
the Segre embedding by tensoring with a generic $\O_H$, where $H$ is a
hyperplane.  (v) follows from (ii)-(iv) and the following general remark:
Let $I\subset S:=k[x_0,\dots,x_n]$ be the ideal of a projective scheme
$X\subset \PP^n:=\Proj S$.  Let
$T:=k[x_0,\dots,x_n,x_{n+1},\dots,x_{n+m}]$ be the flat $S$-algebra
such that $\PP^{n+m}=\Proj T$.  Then the cone over $X$ in $\PP^{n+m}$
has ideal $I\otimes_S T$. Since $T$ is flat over $S$, the resolution
(i.e. the graded betti numbers)
of $I$ is preserved; just exchange $\O_{\PP^{n+m}}$ for $\O_{\PP^n}$.
\end{proof}

For reducible schemes the situation is slightly more difficult to
state, although the geometric picture is clear.  The new subtlety is
the intersections of the irreducible components.  We first state the
precise result describing the possibilities for the union of
four-dimensional irreducible varieties $L$ and $Q$ of degree one and
two:

\begin{prop}
The possible schemes that are a union of $L$ and $Q$ are classified by
the intersection $L\cap Q$.  The intersection is one of the following:

\begin{itemize}
\item[(i)] The empty set (i.e. the disjoint union).
\item[(ii)] One point.
\item[(iii)] Two points.
\item[(iv)] A line.
\item[(v)] A conic.
\item[(vi)] A plane.
\item[(vii)] A two dimensional quadric.
\item[(viii)] A three-dimensional linear space.
\item[(ix)] A three dimensional quadric.
\end{itemize}

Thus the corresponding Hilbert polynomials are (recall that the
Hilbert polynomials of $L$ and $Q$, in our notation, are $P_4$ and
$2P_4-P_3$) 
\begin{itemize}
\item[(i)] $3P_4-P_3$.
\item[(ii)] $3P_4-P_3-P_0$.
\item[(iii)] $3P_4-P_3-2P_0$.
\item[(iv)] $3P_4-P_3-P_1$.
\item[(v)] $3P_4-P_3-2P_1+P_0$.
\item[(vi)] $3P_4-P_3-P_2$.
\item[(vii)] $3P_4-P_3-2P_2+P_1$.
\item[(viii)] $3P_4-2P_3$.
\item[(ix)] $3P_4-3P_3+P_2$.
\end{itemize}
\end{prop}
\begin{proof}
The first part is a simple consequence of the fact that the highest
dimensional linear subvarieties lying on an irreducible quadric hypersurface in
$\PP^5$ have dimension two.  The second part follows from the first
and the equality
\[ \Hilb (L\cup Q)=\Hilb L +\Hilb Q - \Hilb(L\cap Q).\]
\end{proof}

We will need the following proposition:

\begin{prop}
Let $W=X\cup Y \cup Z\subset \PP^N$ be a union of three pairwise different
four-dimensional linear subspaces.  Let the Hilbert polynomial of $W$
be
\[\Hilb W=3P_4-aP_3+\dots\]
Then $a=0,1,2$ or $3$.
\end{prop}

\begin{proof}
Intersect $W$ with a general $\PP^{N-3}$.  This intersection is then a
union of three lines with Hilbert polynomial
\[3P_1-aP_0.\]

The possible unions of three lines are the following:
\begin{itemize}
\item Three skew lines, giving $a=0$.
\item Two lines intersecting in a point, the third not intersecting
  these two, giving $a=1$.
\item Two skew lines, the third intersecting both these two, giving
  $a=2$.
\item A triangle, giving $a=3$.
\item Three lines in the plane, intersecting in one point, giving
  $a=3$.
\item Three non-coplanar lines intersecting in one point, giving $a=2$.
\end{itemize}
This establishes the claim.
\end{proof}

\begin{bem} In fact, we can be more precise.  Through a case-by-case
  analysis it can be seen that the following is the complete list of
  Hilbert polynomials for a union of three four-dimensional linear
  spaces:

\[\begin{array}{lll}
3P_4&3P_4-P_0&3P_4-2P_0\\
3P_4-3P_0&3P_4-P_1&3P_4-P_1-P_0\\
3P_4-2P_1&3P_4-P_1-2P_0&3P_4-2P-1-P_0\\
3P_4-3P_1+P_0&3P_4-3P_1&3P_4-P_2\\
3P_4-P_2-P_0&3P_4-P_2-P_1&3P_4-P_2-2P_0\\
3P_4-P_2-P_1&3P_4-P_2-P_1-P_0&3P_4-2P_2\\
3P_4-P_2-2P_1+P_0&3P_4-P_2-2P_1&3P_4-2P_2-P_1+P_0\\
3P_4-3P_2+P_1&3P_4-3P_2+P_0&3P_4-P_3\\
3P_4-P_3-P_0&3P_4-P_3-2P_0&3P_4-P_3-P_1\\
3P_4-P_3-2P_1+P_0&3P_4-P_3-P_2&3P_4-P_3-2P_2+P_1\\
3P_4-2P_3&3P_4-3P_3+P_2&
\end{array}\]

Since the argument is quite technical, and we only need the weaker
statement in the above proposition, we will not provide the proof for
this statement.
\end{bem}

\begin{prop}
Let $Y$ be the triple structure on $\PP^4$ embedded in $\PP^{39}$ from 
Example \ref{HorrMum}.  Then $Y$ determines a point in a Hilbert
scheme $H$.  No point in $H$ corresponds to a reduced,
equidimensional scheme.
\end{prop}
\begin{proof}
The Hilbert polynomial of the Horrocks-Mumford bundle $\E$ was computed
already in \cite{HM73}.  With our twist and notation, the result reads
\[\Hilb \E =2P_4+5P_3+5P_2-10P_0. \]
\noindent
By Proposition \ref{hilbofmult} the Hilbert polynomial of $Y$ is then
\[\Hilb Y =\Hilb \PP^4+\Hilb \E=3P_4+5P_3+5P_2-10P_0\]
\noindent
which is not among the possibilities for reduced degree three schemes
found above.
\end{proof}

\begin{bem}
One might hope to gain new insight into the validity of Hartshorne's
Conjecture through the reformulation in Theorem \ref{Hartshorneconj}
and the very small number
of degree three varieties.  Our last proposition suggests that it will
be difficult to search for (possible) counter-examples to the
conjecture by this procedure.  At least we cannot expect to find new
bundles by deforming given reduced varieties to triple schemes on
linear subspaces and reversing the above process.
\end{bem}

\begin{bem}  In general, for arbitrary dimension $n$, we can prove by
  the same argument that any reduced, equidimensional scheme of degree
  three has Hilbert polynomial of the form $3P_n-aP_{n-1}+\dots$, where
  $a=0,1,2,3$. 
\end{bem}

\subsection{Complete intersections, projections and multiple structures}

  In this section we consider subvarieties of a projective space
  $\PP^N$.  Our first result provides a partial generalization of
  a criterion from \cite{BE90} from codimension two to higher
  codimension.  Roughly speaking, it says that complete intersections are
  characterized by the structure of the first infinitesimal
  neighbourhood.

\begin{thm}
Let $X\subset \PP^N$ be a smooth variety of codimension $c$.
\begin{itemize}
\item[a)] Let $X$ be a complete intersection, and choose polynomials
  $F_1,\dots ,F_c$ of degrees $f_1,\dots,f_c$ generating the ideal of
  $X$.  For each subset $S\subset \{1,\dots,c\}$ there is a multiple
  Cohen-Macaulay 
  structure $Z_S$, contained in the first infinitesimal neighbourhood, with
  defining short exact sequence
\[0\rightarrow \J_S/\I_X^2\rightarrow
\I_X/\I_X^2\stackrel{F_i}{\rightarrow}\oplus \O_X(-f_i)\rightarrow 0,\]
where $i$ runs through $S$.  $Z_S\subset Z_{S'}$ if and only if
$S\subset S'$.  The multiplicity of $Z_S$ is equal to the cardinality
of $S+1$.
\item[b)] If $X$ has multiple structures $Z_S$ as above (but with
  general line bundles, not necessarily of the form $O_X(-f)$), then the
  normal bundle of $X$ splits.  Especially, if the Picard group of $X$
  is generated by the class of a hyperplane section (e.g. if $\dim X\geq
  \frac{N+2}{2}$), $X$ is a complete intersection.

\item[c)] Assume that the Picard group of $X$ is generated by the
  class of a
  hyperplane section, and that $H^2_{\ast}\I_X=0$.  Then $X$ is a complete
intersection if and only if there is a filtration 
\[X=X_0\subset X_1\subset ... \subset X_{\codim X}=X^{(1)}\]
of the first infinitesimal neighbourhood of $X$ in $\PP$, where each
$X_i$ is Cohen-Macaulay and has multiplicity one more than $X_{i-1}$.
\end{itemize}
\end{thm}
\begin{proof}
The proof of part a): the conormal bundle of $X$ is
isomorphic to $\oplus_i\O_X(-f_i)$, and the morphisms in the statement
are projections onto summands given by choosing the polynomials not in
the ideal of $Z_S$.  The last statements follow easily.\\

For the proof of part b), assume that the codimension $c=3$.  The extension
to general codimension is straightforward but messy.  So we have a
graph of inclusions 
\[\xymatrix{
&Z_1\ar[r]\ar[rd]&Z_{12}\ar[rd]&\\
X=Z_{\emptyset}\ar[ru]\ar[r]\ar[rd]&Z_2\ar[ru]\ar[rd]&Z_{13}\ar[r]&
Z_{123}=X^{(1)} \\
&Z_3\ar[r]\ar[ru]&Z_{23}\ar[ru]&}\]

Now for each pair $i,j$ there is a diagram
\[\xymatrix{
&&&0\ar[d]&\\
&0\ar[d]&&\M\ar[d]&\\
0\ar[r]&\I_{ij}/\I_X^2\ar[r]\ar[d] &\I_X/\I_X^2\ar[r]\ar@{=}[d]&\L_{ij}
\ar[r]\ar[d] & 0\\
0\ar[r]&\I_i/\I_X^2\ar[r]\ar[d] &\I_X/\I_X^2\ar[r]&\L_i\ar[r]\ar[d]&0\\ 
&\M\ar[d]&&0&\\
&0&&&}\]
Reversing the role of $i$ and $j$ gives a similar diagram, but the
the right hand column is turned upside down, showing two
things:  that the $\M$ defined by the diagram is isomorphic to $\L_j$,
and that the right hand sequence splits.  Now a similar argument can be
made for the diagram
\[\xymatrix{
&&&0\ar[d]&\\
&0\ar[d]&&\M'\ar[d]&\\
0\ar[r]&\I_{ijk}/\I_X^2\ar[r]\ar[d] &\I_X/\I_X^2\ar[r]\ar@{=}[d]&\L_{ijk}
\ar[r]\ar[d] & 0\\
0\ar[r]&\I_{ij}/\I_X^2\ar[r]\ar[d] &\I_X/\I_X^2\ar[r]&\L_{ij}\ar[r]\ar[d]&0\\ 
&\M'\ar[d]&&0&\\
&0&&&}\]
showing that $\M'\iso \L_k$ and that $\L_{ijk}$ splits as a sum of the
three line bundles $\L_i,\L_j,\L_k$.  But $\I_{ijk}=\I_X^2$, so 
\[\I_X/\I_X^2\iso \L_{123}\iso \L_1\oplus\L_2\oplus\L_3.\]

Now if the Picard group is generated by the class of a hyperplane section, these
summands are all twists of the structure sheaf of $X$, and the normal
bundle of $X$ can be extended to the corresponding split bundle on
$\PP^N$.\\

The condition $\dim X\geq \frac{N+2}{2}$ is the condition in the
Barth-Lefschetz Theorem, which says, among other things, that $\Pic
X=h\ZZ$, $h$ the class of a hyperplane section.\\

Let us prove part c).  
Suppose we have a filtration as in the statement.  For each $i$ there
is a short exact sequence
\[0\rightarrow \I_{X_i}/\I_X^2\rightarrow \I_X/\I_X^2\rightarrow
\L_i\rightarrow 0\]
 \noindent
where $\L_i$ is locally free of rank $i$.  Furthermore we have
surjections $\L_i\surj \L_j$ whenever $i>j$.  Thus we have a
filtration of $\I_X/\I_X^2$ with successive quotients which are line
bundles on $X$.  Because of the conditions
$\Ext^1(\M,\N)=H^1\O_X(d)=H^2\I_X(d)=0$ for all line bundles $\M,\N$ on $X$
($d$ being some twist depending on $\M$ and $\N$),
and $\I_X/\I_X^2$ actually splits.  We conclude as in part b).

\end{proof}

\begin{bem}
Of course, given a filtration as in part c) one can get the splittings if the
cohomology vanishes in the (finite number of) relevant degrees.
\end{bem}

\begin{bem}
In codimension two, the result has been generalized to singular
varieties, see \cite{FKLpre}.
\end{bem}

Our next results tie in double structures and projections.  Our setup
is the following:  $i:X\hookrightarrow \PP^N$ is a non-singular, connected
variety, and $P\in \PP^N\setminus X$ is the center of a projection
$\pi:\PP^N\setminus \{P\} \rightarrow \PP^{N-1}$.  We use the same
symbol for the restriction to $X$.  $X'\subset \PP^{N-1}$ is the image
of $X$ under this morphism.

\begin{thm}  In the above situation, 
\begin{itemize}
\item[a)] If $P$ is outside $\Sec X$, the secant variety of $X$, then there
  is a double Cohen-Macaulay structure $Y$ on $X$ such that
  its defining sequence is
\[0\rightarrow \I_Y/\I_X^2\rightarrow \I_X/\I_X^2\rightarrow
\O_X(-1)\rightarrow 0.\]
  In addition, 
\[\I_Y/\I_X^2\iso \I_{X'}/\I_{X'}^2.\]
\item[b)] Suppose that $X\subset \PP^3$ is a space curve, and that
  $X'$ only has ordinary nodes as singularities.  Then
  there exists a double Cohen-Macaulay structure $Y$ on $X$
  such that its defining sequence is
\[0\rightarrow \I_Y/\I_X^2\rightarrow \I_X/\I_X^2\rightarrow
\O_X(-1)\rightarrow 0.\]
  In addition
\[\I_Y/\I_X^2\iso \O_X(-\deg X')\otimes \O_X(\DD)\]
 where $\DD$ is the divisor of double points of the morphism
 $\pi:X\rightarrow \PP^2$.
\end{itemize}
\end{thm}

\begin{proof} 
a) There is a well-known exact sequence
\[0\rightarrow \I_{X'}/\I_{X'}^2\rightarrow \I_X/\I_X^2\rightarrow
i^{\ast}\Omega^1_{\PP^N\setminus \{P\}/\PP^{N-1}}\rightarrow 0\]
(see \cite{SGA6}, VIII.1).
We use our standard procedure on this sequence to produce the double
structure $Y$.  The only thing left to observe is that
\[i^{\ast}\Omega^1_{\PP^N\setminus \{P\}/\PP^{N-1}}\iso \O_X(-1)\]
which can be seen, for instance, by pulling back the Euler sequences
of $\PP^N$ and $\PP^{N-1}$ to $X$:
\[\xymatrix{
&&0&0&\\
&&\O_X(1)\ar[u]\ar@{=}[r]&\O_X(1)\ar[u]&\\
0\ar[r]&\O_X
\ar[r]\ar@{=}[d]&\O_X(1)^{\oplus N+1}\ar[r]\ar[u]&i^{\ast}\T_{\PP^N}\ar[r]\ar[u]&0\\ 
0\ar[r]&\O_X\ar[r]&\O_X(1)^{\oplus N}\ar[r]\ar[u]&\pi^{\ast}\T_{\PP^{N-1}}\ar[r]\ar[u]&0\\
&&0\ar[u]&0\ar[u]&}\]
The dual of the sequence to the right is then
\[0\rightarrow \O_X(-1)\rightarrow i^{\ast}\Omega_{\PP^N}\rightarrow
\pi^{\ast}\Omega_{\PP^{N-1}} \rightarrow 0\]
from which we get
\[i^{\ast} \Omega^1_{\PP^N\setminus \{P\}/\PP^{N-1}}\iso \O_X(-1).\]

b)  In this situation there is an exact sequence
\[0\rightarrow \T_X\rightarrow \pi^{\ast}\T_{\PP^2}\rightarrow
\N_{\pi}\rightarrow 0,\]
where the relative normal bundle is
\[\N_{\pi}\iso\pi^{\ast}\O_{\PP^2}(\deg X')\otimes \O_X(-\DD).\]
This can be found in \cite{Ful98},9.3.\\

Consider also the exact sequence
\[0\rightarrow \T_X \rightarrow i^{\ast}\T_{\PP^3}\rightarrow
\N_{X/\PP^3}\rightarrow 0,\]
the normal sequence to the embedding $i:X\hookrightarrow \PP^3$.
Combining these sequences with the diagram from the proof of part a)
(but with the arrows reversed),
we get a diagram
\[\xymatrix{
&&0\ar[d]&0\ar[d]&\\
&&\O_X(1)\ar[d]\ar@{=}[r]&\O_X(1)\ar[d]&\\
0\ar[r]&\T_X
\ar[r]\ar@{=}[d]&i^{\ast}\T_{\PP^3}\ar[r]\ar[d]&\N_X\ar[r]\ar[d]&0\\ 
0\ar[r]&\T_X\ar[r]&\pi^{\ast}\T_{\PP^2}\ar[r]\ar[d]&\N_{\pi}\ar[r]\ar[d]&0\\
&&0&0&}\]
\noindent
Now the dual of the sequence to the right,
\[0\rightarrow \N_{\pi}^{\vee} \rightarrow \N_X^{\vee} \rightarrow
\O_X(-1) \rightarrow 0,\]
is the sequence we seek.
\end{proof}

\begin{bem}
The multiple structures defined in the theorem can be given an
explicit geometric interpretation as follows: for each point $x$ on $X$,
construct the line joining this point to the point $P$.  For a point
on a line, there is only one double structure.  The theorem gives
conditions under which these double points patch together to give a
Cohen-Macaulay double structure on $X$.  We double $X$ in the
direction of $P$.\\

The dual to the morphism $\I_X/\I_X^2\rightarrow \O_X(-1)$,
\[\O_X(1)\hookrightarrow \N_X,\]
can be thought of as giving a normal vector in each point of $X$,
pointing towards $P$.  If $P$ were on a tangent line to a point $x$ on
$X$, this vector would vanish, i.e. the cokernel of
$\O_X(1)\rightarrow \N_X$ would have torsion with support on $x$.
This explains why $P$ has to be outside the tangent variety
of $X$ in the theorem (e.g. in part b) we assumed that $X'$ did not
have any cusps).\\

Also, if $X$ is degenerate and $P$ is outside the linear span of $X$,
this $\O_X(-1)$ is a direct summand of $\I_X/\I_X^2$; compare with the
proof of Proposition \ref{buntkvotient}.
\end{bem}


\section{\label{ikkeI} Examples of non-type I structures}
In this chapter we will find some examples of codimension two multiple
structures
whose $S_1$-filtration contains schemes that are not Cohen-Macaulay.
All our examples will contain a complete intersection of multiplicity
two less, and this will imply that the dimension has to be two.  For a
linear subspace, some of these structures were found in my Master's thesis
(hovedfagsoppgave).  For some classes of reduced support, the examples
presented here give a complete list of structures of this kind.  The
two simplest examples were found by Manolache in \cite{Man92}.\\

More precisely, we will prove the following theorem:

\begin{thm}  Let $X\subset \PP^N$ be a smooth variety of codimension
  two and degree $d$, embedded in a projective space of dimension
  $N\geq 4$.  Assume that $X$ is a complete intersection with ideal
  $I_X=(P,Q)$, and let $Y$ be the multiple structure on $X$ given by
  $I_Y=(P^a,Q^b)$.  Let $Y'$ be a Cohen-Macaulay multiple structure on
  $X$, containing $Y$, with multiplicity two higher than the
  multiplicity of $Y$ (i.e. $ab+2$), such that the $S_1$-filtration of
  $Y'$ contains a non-Cohen-Macaulay scheme $Z$ lying between $Y$ and
  $Y'$.  Then $N=4$.  Possibilities for the ideal of $Y'$ are:

\[\begin{array}{r}  
{\mbox If }a=1,b=1\\
I_{Y'}=(g^2(fP+gQ)-lP^2,fg(fP+gQ)+lPQ,f^2(fP+gQ)-lQ^2, P(fP+gQ),\\
 Q(fP+gQ),P^3,P^2Q,PQ^2,Q^3)\\

{\mbox If }a=1, b>1\\
I_{Y'}=(g(fP+gQ^b)-lPQ,f(fP+gQ^b)+lQ^{b+1},\\
  Q(fP+gQ^b),P^2,PQ^2,Q^{b+2})\\

{\mbox If }a>1,b>1\\
I_{Y'}=(rg(fP^a+gQ^b)-lP^{a+1},
rf(fP^a+gQ^b)+lPQ^b,sg(fP^a+gQ^b)-lP^aQ,\\
sf(fP^a+gQ^b)+lQ^{b+1},sP^{a+1}-rP^aQ, sPQ^b-rQ^{b+1},\\
P(fP^a+gQ^b),Q(fP^a+gQ^b),P^{a+2},P^2Q^b,P^{a+1}Q,PQ^{b+1},
P^aQ^2,Q^{b+2})
\end{array}\]

Furthermore, if $\Pic X\iso \ZZ$, generated by a hyperplane section,
then the above list is complete.
\end{thm}

\label{setup}
Let $X\subset \PP^N$ be a smooth variety of codimension two and degree
$d$, $N\geq 4$.  We
assume that $\Pic(X)=h\ZZ$, where $h$ is the class of a hyperplane
section, which holds for instance
if $X$ has dimension at least four, or if $X$ is a complete
intersection of dimension three.  We assume
furthermore that $H^1_{\ast}(\I_X)=0$ (which also holds if $X$ is a complete
intersection, or more generally if the length of a minimal resolution
of $X$ is less than the maximal possible length).  Together, these two
assumptions ensure that all
$\O_X$-Module maps between line bundles lift to morphisms of line
bundles on $\PP^N$, and are therefore represented by polynomials.
Note also that if we remove the condition that the Picard group is
generated by a hyperplane section, or the condition $H^1_{\ast}(\I_X)=0$, the constructions still work,
although there may be other structures arising from constructions
using line bundles or maps that
cannot be extended to $\PP^N$.  Let $Y$
be a complete intersection with $Y_{\red}=X$, $I_Y=(p,q)\subset
S=k[z_0,...,z_N]$.  We want to understand short exact sequences

\[0\rightarrow \I_Z/\I_Y\I_X\rightarrow \I_Y/I_Y\I_X \rightarrow \L
\rightarrow 0\]
\noindent
where $Z$ is an $(S_1)$-scheme, but not Cohen-Macaulay.  Since $Y$ is a
complete intersection, $\I_Y/\I_Y\I_X$ is isomorphic to
$\O_X(-\alpha)\oplus\O_X(-\beta)$, where $\alpha,\beta$ are the degrees of $p,q$
respectively.  $\L$ is torsion free of rank 1, so $\I_Z/\I_Y\I_X$ is
reflexive and thus, being of rank one, a line bundle.  If the sequence
splits, $\L$ has to be a line bundle and we get a Cohen-Macaulay
scheme. 
  
\begin{bem}
In this opening
section the polynomials $p$ and $q$ are just required to be
coprime, so that the scheme they define has codimension two.  Later on
we will also require the reduced subvariety to be a complete
intersection, and $p$ and $q$ will be powers of generators of the
ideal of $X$.  This is to facilitate the computations; it seems harder
to do this in a more general setting.
\end{bem}

\begin{prop}
In order to get a structure $Z$ of degree $\alpha \beta+d$ (or equivalently, of
multiplicity $\alpha \beta/d+1$) that contains $Y$, which is $(S_1)$ but not
Cohen-Macaulay, it is necessary and sufficient to require that 
\[\I_Z/\I_Y\I_X\iso \O_X(-c),\mbox{  }c>\alpha,\beta.\]
\end{prop}

\begin{proof}
We already know that $\I_Z/\I_Y\I_X$ has to be a line bundle, so we
just have to consider the twist $-c$.  A necessary condition for $Z$
is that the homological singularity set has codimension at least
two.  From the existence of a 
non-zero map $\O_X(-c)\iso\I_Z/\I_Y\I_X\rightarrow
\I_Y/I_Y\I_X\iso\O_X(-\alpha)\oplus\O_X(-\beta)$ we conclude that $c\geq$ one
of $\alpha,\beta$.  Suppose $\alpha\geq \beta$.  \\

If $\alpha>c\geq \beta$ the first component of the map
$\O_X(-c)\rightarrow\O_X(-\alpha)\oplus\O_X(-\beta)$  will be zero.  This
violates the condition of codimension two for the homological
singularity set.  Likewise, if $c=\alpha$ the first component of our map
will be invertible (or zero), again violating this condition.  Thus
$c>\alpha,\beta$ is necessary.\\

For the sufficiency, note that the map
$\O_X(-c)\rightarrow\O_X(-\alpha)\oplus\O_X(-\beta)$ will be given by two forms
$f,g$ of degrees $c-\alpha,c-\beta$ respectively, such that $p,q,f,g$ is a
regular sequence.  Obviously $Z$ will satisfy
the Cohen-Macaulay condition in 
each point outside $V_+(f,g)$, so we need only show that the
$(S_1)$-condition is satisfied in the points of $V_+(f,g)$.  In these
points $P$ we have $\dim \O_{Z,P}\geq 2$.  By definition, we just have
to show that $\depth \O_{Z,P}\geq 1$.  But to
show this it suffices to present a non-zero divisor modulo $I_Z$, for
instance $f$ or $g$.
\end{proof}

If the ideal of $X$ is $I_X=(P,Q)$ then, in terms of $f,g$ the ideal
of $Z$ is $I_Z=(fp+gq,pP,pQ,qP,qQ)$.  

\begin{prop}
\label{Sisteerlinjebunt}
If the last term in the $S_1$-filtration, $Z$, of a Cohen-Macaulay scheme
$Y$, of multiplicity $\nu$, is of multiplicity
$\nu-1$, then the last $\L_j$ is a line bundle.
\end{prop}
\begin{proof}
Let $\L=\L_j$.  There is a short exact sequence of
  $\O_{\PP^N}$-Modules  

\[0\rightarrow i_{\ast}\L\rightarrow \O_Y \rightarrow \O_Z \rightarrow
0\]
\noindent
where $i:X\hookrightarrow \PP^N$ is the inclusion.  By the results on
completion, we need to show that, for each $x\in X$, the
$A:=\widehat{\O_{X,x}}$-module $L:=\widehat{\L_x}$ is free.
Now this module sits in a short exact sequence

\[0\rightarrow L \rightarrow  \widehat{\O_{Y,x}}\rightarrow
\widehat{\O_{Z,x}} \rightarrow 0\]
\noindent
of $A$-modules.  In this exact sequence the module in the middle is
free and the module on the right is torsion free.  Thus $L$ is
reflexive, and, being of rank one, also free.  Since $x\in X$ was
arbitrary, the sheaf $\L$ is a line bundle. 
\end{proof}

For the explicit calculations that follow, we will specialize to a
complete intersection with ideal $I_Y=(p,q)=(P^a,Q^b)$, where
$I_X=(P,Q)$ is the ideal of the support $X$ (also a complete
intersection). Thus 
$I_Z=(fP^a+gQ^b,P^{a+1},PQ^b,P^aQ,Q^{b+1})$, 
$I_ZI_X=(P(fP^a+gQ^b),Q(fP^a+gQ^b),P^{a+2},P^2Q^b,P^{a+1}Q,PQ^{b+1},
P^aQ^2,Q^{b+2})$. Throughout, $\O=\O_X$.\\

\begin{bem} Note that if $a=1$ or $b=1$ there are redundant generators
  in the above expression for $I_Z$.  Therefore we must split up our
  analysis accordingly.  However, it turns out that the required
  calculations are very similar for varying cases.  Therefore we give
  the complete proof only for the case that both $a$ and $b$ are
  greater than two.  For further details, see \cite{Vat01}.
\end{bem}

\begin{prop} Assume both $a\geq 2$ and $b\geq 2$.  The minimal
  resolution of the $\O$-Module $\I_Z/\I_Z\I_X$ is 

\[ 0\rightarrow \begin{array}{c} \O(-c-\deg P) \\ \oplus \\
  \O(-c-\deg Q) \end{array} \rightarrow
\begin{array}{c} \O(-c)\\ \oplus\\ \O(-\alpha-\deg P)\\ \oplus\\
  \O(-\beta -\deg P) \\ \oplus \\ \O(-\alpha -\deg Q) \\ \oplus \\
  \O(-\beta-\deg Q)\end{array} \rightarrow \I_Z/\I_Z\I_X \rightarrow 0
  \] 
\noindent
where $\alpha=\deg(P^a)=a\deg P$,$\beta =b\deg Q$ and
$c=\deg(fP^a+gQ^b)$.  The first map is given, on the level of graded
modules, by the matrix
\[ \Phi=\left( \begin{array}{cc}
0&0\\
f&0\\
g&0\\
0&f\\
0&g
\end{array} \right) \]
\noindent
and the second map is given by the generators of $I_Z$
\[\Psi=\left( \begin{array}{ccccc}
fP^a+gQ^b & P^{a+1} & PQ^b & P^aQ & Q^{b+1}\end{array}\right). \]
\end{prop}

\begin{proof}  The injectivity of the first map and the surjectivity
  of the second is obvious.  First we show that $\Psi \circ \Phi=0$
  in $\O$-Mod. 

\[\begin{array}{l}
\Psi \circ \Phi=\left( \begin{array}{ccccc}
fP^a+gQ^b & P^{a+1} & PQ^b & P^aQ & Q^{b+1}\end{array}\right) \left( \begin{array}{cc}
0&0\\
f&0\\
g&0\\
0&f\\
0&g
\end{array} \right)= \\
=\left( \begin{array}{cc} fP^{a+1}+gPQ^b & fP^aQ+gQ^{b+1}\end{array}
\right). \end{array}\]
Now both these polynomials lie in $I_ZI_X$, so the composition is
zero.\\

Then we show that $\ker \Psi\subset \im \Phi$.  So suppose
$\Phi(A,B,C,D,E)=0$, that is 
\[A(fP^a+gQ^b)+BP^{a+1}+CPQ^b+DP^aQ+EQ^{b+1}=0\]
or, if we consider preimages of $A,B,C,D,E$ (which we denote by the
same letters) 
\[A(fP^a+gQ^b)+BP^{a+1}+CPQ^b+DP^aQ+EQ^{b+1}\in I_ZI_X.\]
Since we are free to adjust these polynomials with elements from
$I_X$, we can assume that each one is either outside $I_X$ or else is
zero.  Let us first see that $A=0$: since $f,A\notin I_X$ the term
$AfP^a$ is indivisible by both $P^{a+1}$ and $Q$. But all other terms
that involve $P^a$ also involve one more $P$ or $Q$.  Thus $A=0$.\\

Consider $B$ and $C$: we want to show that $f|B,g|C$ and
$\frac{B}{f}=\frac{C}{g}$ (the analogous statement for $D$ and $E$ is
proven similarly). Since 
$B\notin I_X$ the term $BP^{a+1}$ is indivisible by $P^{a+2}$ and $Q$.
Since $C\notin I_X$ the term $CPQ^b$ is indivisible by $P^2$ and
$Q^{b+1}$.  Since all generators for $I_ZI_X$, apart from  
$P(fP^a+gQ^b)$, involve higher powers of $P$ and $Q$, we will have to
use this relation.  Thus there is some polynomial $h$ such that
$BP^{a+1}+CPQ^b=h(fP^{a+1}+gPQ^b)$ or

\[BP^a+CQ^b=h(fP^a+gQ^b).\]
\noindent
Since $P,Q$ is a regular sequence we can split this equation in two:
\[BP^a=hfP^a\mbox{ and }CQ^b=hgQ^b.\]
Finally,
\[\frac{B}{f}=h=\frac{C}{g}\]
and the proof is complete.
\end{proof}

Note that we used the condition $a\geq 2$ in the last part of the
proof.  The condition $b\geq 2$ is used similarly for the $D,E$ case.
\\

The next step is to find the sheaves $\L$ such that there is a
surjective map $\I_Z/\I_Z\I_X\rightarrow \L$ whose kernel gives the
Ideal of a Cohen-Macaulay scheme.  By Proposition \ref{Sisteerlinjebunt}
$\L$ is a line bundle.  We still assume $a,b\geq 2$.

\begin{prop} Giving a surjection 
\[\I_Z/\I_Z\I_X\rightarrow \L\]
to an invertible sheaf $\L$ is equivalent to giving a surjection  

\[\begin{array}{c} \O(-c)\\ \oplus\\ \O(-\alpha-\deg P)\\ \oplus\\
  \O(-\beta -\deg P) \\ \oplus \\ \O(-\alpha -\deg Q) \\ \oplus \\
  \O(-\beta-\deg Q)\end{array} \rightarrow \L \]
\noindent
whose matrix on the level of graded modules is
\[\Lambda= \left( \begin{array}{ccccc} l&rg&-rf&sg&-sf\end{array}\right)\] 
where $r,s,l$ are homogeneous polynomials of degree $\deg r$,$\deg
s=\deg r+\deg Q-\deg P$ and $\deg l=d+c=2c-\alpha-\beta-\deg P+\deg r$
that satisfy
\[V_+(l)\cap V_+(rg)\cap V_+(rf)\cap V_+(sg)\cap V_+(sf)=\emptyset.\]
If such a surjection exists, then $\dim X=2$ and
$\L\iso\O(c-\alpha-\beta-\deg P+\deg r)$.
\end{prop}
\begin{proof}
Since the displayed rank five Module surjects onto $\I_Z/\I_Z\I_X$, we
need a map from this Module that annihilates the kernel, that is 
\[\Lambda\circ \Phi=0.\]
Since we are interested in $\L$ a line bundle, $\Lambda$ will be a
$1\times 5$ matrix.  The map is given by five polynomials, and the
surjectivity will then amount to saying that those five polynomials
have no common zeroes.  Let $\Lambda=\left( \begin{array}{ccccc}
    l&m&n&o&p\end{array}\right) $.  The condition $\Lambda \Phi=0$
translates to the two equations
\[\begin{array}{ccc}
fm+gn&=&0\\
fo+gp&=&0. \end{array} \]
\noindent
This gives the following conditions:
\[\begin{array}{ccc}
l&&\mbox{no condition}\\
m&=&rg\\
n&=&-rf\\
o&=&sg\\
p&=&-sf \end{array}\]
\noindent
The surjectivity condition then amounts to
\[V_+(l)\cap V_+(rg)\cap V_+(rf)\cap V_+(sg)\cap V_+(sf)=\emptyset\]
or
\[\begin{array}{ccc}
V_+(l)\cap V_+(r)\cap V_+(s)&=&\emptyset\\
V_+(l)\cap V_+(f)\cap V_+(g)&=&\emptyset. \end{array}\]
\noindent
If $\dim X>2$ then this is impossible, whereas if $\dim X=2$ then this
will hold for generic choices of polynomials.\\

The only thing left to compute is the degree of $\L$ and the degrees
of the polynomials $r,s$ and $l$; let $\L=\O(d)$.
From the second term of $\Lambda$ we find
\[d= -\alpha-\deg P +\deg r+\deg g=-\alpha-\deg P +\deg r +c-\beta\]
which gives the formula for $d$ in the statement.
From the fourth term of $\Lambda$ we find
\[d=-\alpha -\deg Q +\deg s+\deg g\]
which, compared with the first expression for $d$, gives
\[\deg s=\deg r +\deg Q -\deg P.\]
Finally, from the first term of $\Lambda$ we get the equation
\[d=-c+\deg l\]
which, compared with the first expression for $d$, gives
\[\deg l=d+c=2c+\deg r-\alpha-\beta-\deg P.\]
\end{proof}

The final step is to find the kernel of $\Lambda$, and interpret this
as a statement about equations.

\begin{prop}
The equations we find from the kernel of $\Lambda$, together with the
equations generating $\I_Z\I_X$ give us the ideal of a
Cohen-Macaulay scheme $Y'$:
\[\begin{array}{r}I_{Y'}=(rg(fP^a+gQ^b)-lP^{a+1},
rf(fP^a+gQ^b)+lPQ^b,sg(fP^a+gQ^b)-lP^aQ,\\
sf(fP^a+gQ^b)+lQ^{b+1},sP^{a+1}-rP^aQ, sPQ^b-rQ^{b+1},\\
P(fP^a+gQ^b),Q(fP^a+gQ^b),P^{a+2},P^2Q^b,P^{a+1}Q,PQ^{b+1},
P^aQ^2,Q^{b+2}) \end{array}\]
\end{prop}
\begin{proof}
First we compute the kernel of $\Lambda$; so suppose
$(A,B,C,D,E)\mapsto 0$ under $\Lambda$.  Note that, because of $l,r,s$
and $l,f,g$ being regular sequences, it is enough to consider only
quintuples $(A,B,C,D,E)$ where two of the terms are non-zero.  We get
the following list of $(A,B,C,D,E)$ generating the kernel:
\[\begin{array}{c}
(rg,-l,0,0,0)\\
(rf,0,l,0,0)\\
(sg,0,0,-l,0)\\
(sf,0,0,0,l)\\
(0,f,g,0,0)\\
(0,s,0,-r,0)\\
(0,sf,0,0,rg)\\
(0,0,sg,rf,0)\\
(0,0,s,0,-r)\\
(0,0,0,f,g) \end{array}\]
\noindent
Note that $(0,sf,0,0,rg)=f(0,s,0,-r,0)+r(0,0,0,f,g)$ and
$(0,0,sg,rf,0)=g(0,0,s,0,-r)+r(0,0,0,f,g)$, so these two quintuples
are superfluous. Note also, though we don't need this fact, that the first four
equations generate the kernel over the quotient field.  For instance,
$(0,0,s,0,-r)= \frac{s}{l}(rf,0,l,0,0)-\frac{r}{l}(sf,0,0,0,l)$.  \\

Combining the eight necessary equations from the above list with the
map $\Psi=\left( \begin{array}{ccccc}
fP^a+gQ^b & P^{a+1} & PQ^b & P^aQ & Q^{b+1}\end{array}\right)$, and
adding the generators for $I_ZI_X$, we get
\[\begin{array}{r}I_{Y'}=(rg(fP^a+gQ^b)-lP^{a+1},
rf(fP^a+gQ^b)+lPQ^b,sg(fP^a+gQ^b)-lP^aQ,\\
sf(fP^a+gQ^b)+lQ^{b+1},fP^{a+1}+gPQ^b, sP^{a+1}-rP^aQ, sPQ^b-rQ^{b+1},\\
fP^aQ+gQ^{b+1},P(fP^a+gQ^b),Q(fP^a+gQ^b),P^{a+2},P^2Q^b,P^{a+1}Q,PQ^{b+1},\\
P^aQ^2,Q^{b+2}). \end{array}\]
Here we have two repetitions coming from the fact that
$fP^{a+1}+gPQ^b=P(fP^a+gQ^b)$ already is in $I_ZI_X$, and similarly
for $fP^aQ+gQ^{b+1}=Q(fP^a+gQ^b)$.  Removing these two equations gives
us the theorem.
\end{proof}

\begin{bem} This completes the proof of the main theorem of this
  Chapter.  For $X=\PP^2$, $a=1$ and $b=1,2$, these structures where
  discovered by Manolache \cite{Man92}.  Complete calculations for the
  two cases left out, for $X=\PP^2$, can be found in \cite{Vat98}.
\end{bem}

\newpage
\bibliography{ref}
\bibliographystyle{plain}
\end{document}